\documentclass[11pt,a4paper]{article}
\usepackage{geometry}
\usepackage[british]{babel} 
\usepackage[T1]{fontenc}
\usepackage[utf8]{inputenc}
\usepackage[pdftex]{graphicx} 
\usepackage{color}
\usepackage{url}
\usepackage{marvosym}
\usepackage{fourier,charter}
\usepackage{varioref}\labelformat{equation}{(#1)}
\usepackage{framed} \definecolor{shadecolor}{gray}{0.90}
\usepackage[round]{natbib}
\usepackage{enumerate,mathrsfs}
\usepackage[modulo,mathlines]{lineno}

\newcommand{\GH}{\mathrm{GH}}
\title{\bf A formulation for continuous mixtures \\
   of multivariate normal distributions}
\author{
  {\Large Reinaldo B. Arellano-Valle}\\
   Departamento de Estadística\\
  Pontificia Universidad Católica de Chile\\
  Chile
  \and
  {\Large Adelchi Azzalini} \\  \large
  Dipartimento di Scienze Statistiche \\
  Università di Padova\\
  Italia
 }
\date{\today}
\let\phi=\varphi
\newcommand{\GMN}{\mathrm{GMN}}

\usepackage{amsfonts}

\newcommand{\equald}{\stackrel{d}{=}}
\newcommand{\E}[2][]{
   \ensuremath{\mathbb{E}_{#1}\!\left\{\displaystyle{#2}\right\}}}
\newcommand{\half}{\mbox{$\textstyle \frac{1}{2}$}}   
\newcommand{\pr}[2][]{
   \ensuremath{\mathbb{P}_{#1}\!\left\{\displaystyle{#2}\right\}}}
\newcommand{\Real}{\mathbb{R}}
\newcommand{\T}{^{\top}}
\newcommand{\inv}{^{-1}}
\newcommand{\var}[2][]{
   \ensuremath{\textrm{var}_{#1}\!\left\{\displaystyle{#2}\right\}}}
\renewcommand{\pr}[2][]{
   \ensuremath{\mathbb{P}_{#1}\!\left\{\displaystyle{#2}\right\}}} 
\newcommand{\N}{\mathrm{N}}     
\newcommand{\SN}{\mathrm{SN}}

\renewcommand{\d}{\,\mathrm{d}}
\newcommand{\dfrac}[2]{\displaystyle{\frac{#1}{#2}}}
\newtheorem{theorem}{Theorem}

\newtheorem{proposition}[theorem]{Proposition}
\newcommand{\tr}{\mbox{\rm tr}}

\begin{document}
\modulolinenumbers[2]
\maketitle
\begin{abstract}
Several formulations have long existed in the literature in the form of 
continuous mixtures of  normal variables where a mixing variable operates 
on the mean  or  on the variance or on both the mean and the variance 
of a multivariate normal variable, by changing the nature of these basic 
constituents from constants to random quantities.
More recently, other mixture-type constructions have been introduced,
where the core random component, on which the mixing operation operates,
is not necessarily normal.
The main aim of the present work is to show that many existing
constructions can be encompassed by a formulation where normal
variables are mixed using two univariate random variables.
For this  formulation, we derive various  general properties.
Within the proposed framework, it is also simpler to formulate new
proposals of parametric families and we provide a few such instances.
At the same time, the exposition provides a review of the theme of
normal mixtures.
\end{abstract}

Key-words: location-scale mixtures, mixtures of normal distribution.
\par

\clearpage
\section{Continuous mixtures of normal distributions}  \label{s:normal-mix}

In the last few decades, a number of formulations have been put forward,
in the context of distribution theory, where a multivariate normal variable
represents the basic constituent but with the superposition of another
random component,  either in the sense that the normal mean value
or the variance matrix or both these components are subject
to the effect of another random variable of continuous type.
We shall refer to these constructions as `mixtures of normal variables';
the matching phrase `mixtures di normal distributions' will also be used.

To better focus ideas, recall a few classical instances of the delineated scheme.
Presumably, the best-known such formulation is represented by scale mixtures of
normal variables, which can be expressed as
\begin{equation}
   Y = \xi + V^{1/2} \, X
   \label{e:var-mix}
\end{equation}
where $X\sim\N_d(0,\Sigma)$, $V$ is an independent random variable on $\Real^+$,
and $\xi\in\Real^d$ is a vector of constants.
Scale mixtures  \ref{e:var-mix} provide a stochastic representation of
a wide subset of the class of elliptically contoured distributions,
often called briefly elliptical distributions.
For a standard account of elliptical distributions,
see for instance \citet{fang:kotz:ng:1990};
specifically, their Section~2.6 examines the connection with scale mixtures
of normal variables.
A very importance instance occurs when $1/V \sim \chi^2_\nu/\nu$,
which leads to the multivariate Student's $t$ distribution.

Another very important construction is the normal variance-mean mixture
proposed by \citet{barndorff-nielsen:1977prs,barndorff-nielsen:1978sjs}
and extensively developed by subsequent literature, namely
\begin{equation}
   Y = \xi + V\, \gamma + V^{1/2} \, X
   \label{e:var-mean-mix}
\end{equation}
where $\gamma\in\Real^d$ is a vector of constants and $V$ is assumed
to have a generalized inverse Gaussian (GIG) distribution.
In this case $Y$ turns out to have a generalized hyperbolic (GH) distribution,
which will  recur later in the paper.

Besides \ref{e:var-mix} and \ref{e:var-mean-mix}, there exists a multitude
of other constructions which belong to the  idea of normal mixtures
delineated in the opening paragraph.
Many of these formulations will be recalled in the subsequent pages,
to illustrate the main target of the present contribution,
which is to present a general formulation for normal mixtures.
Our proposal involves an additional random component, denoted $U$,
and the effect of $U$ and $V$ is regulated by two functions,
non-linear in general.
As we shall see, this construction  encompasses a large number of
existing constructions in a unifying scheme, for which we develop 
various general properties.

The role of this activity is to highlight the relative connections of the
individual constructions, with an improved understanding of their nature.
As a side-effect, the presentation of the individual formulations plays
also the role of a review of this stream of literature.
Finally, the proposed formulation can facilitate the conception of
additional proposals with specific aims.
The emphasis is primarily on the multivariate context.

Since it moves a step towards generality,  we mention beforehand
the formulation of \cite{tjet:sene:2006} where  $V$ and $V^{1/2}$  in
\ref{e:var-mean-mix} are replaced  by two linear functions of them,
which allows to incorporate a number of existing families.
Their construction is, however, entirely within the univariate
domain. A number of multivariate constructions aiming at some level 
of generality do exist, and will examined in the course of the discussion.

In the next section, our proposed general scheme is introduced,
followed by the derivation of a number of general properties.
The subsequent sections show how to frame a large number of existing
constructions  within the proposed scheme. 
In the final section, we indicate some directions for even more
general constructions.

\section{Generalized mixtures of normal distributions}   \label{s:gmn}


\subsection{Notation and other formal preliminaries}    \label{s:notation}

As  already effectively employed,
the notation  $W\sim\N_d(\mu, \Sigma)$ indicates that $W$ is a $d$-dimensional
normal random variable with mean vector $\mu$ and variance matrix $\Sigma$.
The density function and the distribution function of $W$ at $x\in\Real^d$
are denoted by $\phi_d(x;\mu,\Sigma)$ and $\Phi_d(x;\mu,\Sigma)$.
Hence, specifically, we have
\[
  \phi_d(x;\mu,\Sigma)= \frac{1}{\det(2\pi\:\Sigma)^{1/2}}\:
     \exp\left\{-\half (x-\mu)\T\Sigma\inv (x-\mu)\right\}
\]
if $\Sigma>0$.
When $d=1$, we drop the subscript $d$. When $d=1$ and, in addition, $\mu=0$ and
$\Sigma=1$, we use the simplified notation $\phi(\cdot)$ and $\Phi(\cdot)$
for the density function and the distribution function.

A quantity arising in connection with the multivariate normal distribution,
but not only there, is the Mahalanobis distance, defined
(in the non-singular case) as
\begin{equation}
  \|x\|_\Sigma = \left(x\T \Sigma\inv x\right)^{1/2}\,   \label{e:M-dist}
\end{equation}
which is written in the simplified form $\|x\|$ when $\Sigma$ is
the identity matrix.

A function which will appear in various expressions is the inverse Mills ratio
\begin{equation}
   \zeta(t) =\frac{\phi(t)}{\Phi(t)}, \qquad t\in\Real\,.
   \label{e:zeta}
\end{equation} 

A positive continuous random variable $V$ has a GIG distribution 
if its density function can be written as
\begin{eqnarray}\label{e:pdf-gig}
  g(v;\lambda,\chi,\psi)=
    \frac{\left(\sqrt{\psi/\chi}\,\right)^{\lambda}}%
          {2K_\lambda\left(\sqrt{\chi\psi}\,\right)}
    \:v^{\lambda-1}\:\exp\left(-\frac{1}{2}(\chi\:v\inv+\psi\:v)\right),
    \hspace{3em} v>0,
\end{eqnarray}
where $\lambda\in\Real$, $\psi>0$, $\chi>0$ and $K_\lambda$ denotes
the modified Bessel function of the third kind.
In this case, we write $V\sim \N^{-}(\lambda,\chi,\psi)$.
The numerous properties of the GIG distribution  and interconnections
with other parametric families are reviewed by \cite{jorgensen:1982}.
We recall two basic properties: both the distribution of $1/V$ and
of $c\,V$ for $c>0$ are still of GIG type.
A fact to be used later is that the Gamma distribution is obtained
when $\lambda>0$ and $\chi\to0$.

A result in matrix theory which will be used repeatedly is the
Sherman-Morrison formula for matrix inversion, which states
\begin{equation}\label{e:SM}
   (A+ b\,d\T)\inv = A\inv - \frac{1}{1+d\T A\inv b} A\inv b d\T A\inv
\end{equation}
provided that the square matrix $A$ and the vectors $b, d$ have
conformable dimensions,  and the inverse matrices  exist.

\subsection{Definition and basic facts}    \label{s:basic-facts}

Consider a $d$-dimensional random variable $X\sim\N_d(0,\Sigma)$ and
univariate random variables $U$ and $V$ with joint distribution function
$G(u,v)$, such that $(X, U, V)$  are mutually independent;
hence $G$ can be factorized as $G(u,v)=G_U(u)\,G_V(v)$.
We assume $\Sigma>0$ to avoid technical complications and concentrate
on the constructive process.
These definitions and assumptions will be retained for the rest
of the paper.

Given any real-valued function $r(u,v)$, a positive-valued function $s(u,v)$,
and  vectors  $\xi$ and $\gamma$ in $\Real^d$,
we shall refer to
\begin{eqnarray}
  Y &=& \xi + r(U,V)\,\gamma + s(U,V)\,X   \label{e:gmn}\\
    &=&  \xi + R\,\gamma + S\, X  \label{e:gmn-RS}
\end{eqnarray}
as a generalized mixture of normal (GMN) variables; we have written
$R=r(U,V)$ and $S=s(U,V)$ with independence of $(R,S)$ from $X$.
Denote by $H$ the joint distribution function of $(R,S)$ implied by $r,s,G$.
The distribution of $Y$ is identified by the notation
$Y\sim \GMN_d(\xi, \Sigma, \gamma,  H)$.

For certain purposes, it is useful to think of $Y$ as generated by the
hierarchical construction
\begin{equation}
\begin{array}{rcl}
   (Y|U=u, V=v) &\sim& \N_d(\xi  + r(u,v)\,\gamma,\: s(u,v)^2\,\Sigma)\,,\\
   (U,V) &\sim& G_U\times G_V \,.
\end{array}
\label{e:gmn-hierar}
\end{equation}
For instance, this representation is convenient for computing the mean vector
and the variance matrix as
\begin{eqnarray}
 \E{Y} &=&\E{\E{Y|U,V}}       \nonumber \\
    &=&\E{\xi+r(U,V)\gamma} \nonumber \\
    &=&\xi+\E{R}\gamma   \label{e:E(Y)}
\end{eqnarray}
provided $\E{R}$ exists, and
\begin{eqnarray}
  \var{Y} &= &\var{\E{Y|U,V}}+\E{\var{Y|U,V}}    \nonumber \\
       &=& \var{\xi+r(U,V)\gamma}+\E{s(U,V)^2\,\Sigma}\nonumber  \\
       &=& \var{R}\,\gamma\gamma\T+\E{S^2}\,\Sigma   \label{e:var(Y)}
\end{eqnarray}
provided $\var{R}$ and $\E{S^2}>0$ exist.
Another use of representation \ref{e:gmn-hierar} is to facilitate the development
of some EM-type algorithm for parameter estimation.

Similarly, by a conditioning argument, it is simple to see that the characteristic
function of $Y$ is
\[
  c(t) = \exp(i t\T\xi)\: \E{c_\N\left(t; r(U,V)\gamma, s(U,V)^2\,\Sigma\right)} \,,
  \qquad\quad t\in\Real,
\]
where $c_\N(t;\mu, \Sigma)$ denotes the characteristic function of
a $\N(\mu, \Sigma)$ variable.
Also, the distribution function of $Y$ is
\begin{equation}
  F(y) = \E{\Phi_d(y; \xi + r(U,V)\gamma, s(U,V)^2\,\Sigma)}\,.
  \label{e:cdf}
\end{equation}

Consider the density function of $Y$, $f(y)$, in the case that $S=s(U,V)$ is a non-null constant.
From \ref{e:gmn-hierar} it follows that
\[
   f(y) = \E[G]{\phi_d(y; \xi + r(U,V)\gamma, s(U,V)^2\,\Sigma)}
        = \E[H]{\varphi_d(y; \xi + R\gamma, S^2\,\Sigma)}
\]
where the first expected value is taken with respect to the distribution $G$, the second one
with respect to $H$.
Assume further that the distribution $H$ of $(R,S)$ is absolutely continuous with 
density function $h(rs,)$, and note that the transformation from $(R,S,X)$ to $(R,S,Y)$ 
is invertible, so that a standard computation
for densities of transformed variables yields, in an obvious notation,
\begin{eqnarray*}
   f_{R,S,Y}(r,s,y) &=& s^{-d}f_{R,S,X}(r,s,s^{-1}(y-\xi-r\gamma)) \\
        &=& h(r,s)\:s^{-d}\:\phi_d(s^{-1}(y-\xi-r\gamma);0,\Sigma) \\
        &=& h(r,s)\:\phi_d(y; \xi+r\gamma, s^2\,\Sigma)
\end{eqnarray*}
taking into account the independence of $(R,S)$ and $X$.
Hence we arrive at
\begin{equation}
   f(y)=\int_{\Real\times\Real_+}\phi_d(y;\xi+r\gamma, s^2\,\Sigma) \d{H(r,s)} \,.
   \label{e:pdf-gmn}
\end{equation}
An alternative route to obtain this expression would be via differentiation of
the distribution function \ref{e:cdf} with exchange of the integration and
differentiation signs.

For statistical work, it is often useful to consider constructions of type \ref{e:gmn}
where the distributions of $U$ and $V$ belong to some parametric family.
In these cases, care must be taken to avoid overparameterization.
Given the enormous variety of  specific instances embraced by \ref{e:gmn},
it seems difficult to establish general suitable condition, and we shall 
then  discuss this issue within specific families or classes of distributions.

In the above passage, as well as in the rest of the paper, the term `family' refers
to the set of distributions obtained by a given specification of the
variables  $(X, U, V)$ when  their parameters vary  in some admissible space,
while keeping  the other ingredients fixed.
Broader sets, generated for instance when the distributions of $U$ and $V$ 
vary across various parametric families,  constitute  `classes'.
 
A clarification is due about the use of the notation in 
\ref{e:gmn}--\ref{e:gmn-RS} and some derived expressions to be presented later on.
When we shall examine a certain family belonging to the general construction, 
that notation will translate into a certain parameterization, 
which often is \emph{not} the most appropriate for inferential 
or for interpretative purposes, and its use here must \emph{not} be 
intended  as a recommendation for general usage.
This scheme is adopted merely for uniformity and simplicity  of 
treatment in the present investigation.


\subsection{Affine transformations and other distributional properties}
\label{s:distrib-properties}

For the random variable $Y$ introduced by \ref{e:gmn}-\ref{e:gmn-RS},
consider an affine transformation $W=b+ B\T\,Y$, for a $q$-dimensional
vector $b$ and a full-rank matrix $B$ of dimension $d\times q$,
with $q\le d$; denote these assumptions as `the b-B conditions'.
It is immediate that
\[
  W=b+ B\T\,Y = b+ B\T\xi + r(U,V)B\T\gamma + s(U,V)\, B\T X
\]
is still of type \ref{e:gmn}--\ref{e:gmn-RS} with the same mixing
variables $(R,S)$ and  modified numerical parameters. We have then
reached the following conclusion.

\begin{proposition} \label{th:affine}
If $Y\sim\GMN_d(\xi, \Sigma, \gamma, H)$ and $b, B$ satisfy  the b-B conditions
introduced above, it follows that
\begin{equation}
  b+B\T Y \sim \GMN_q(b+ B\T\xi, B\T\Sigma B, B\T\gamma, H)
  \label{e:affine}
\end{equation}
is still a member of the GMN class, with the same mixing distribution of $Y$.
\end{proposition}

Partition now $Y$ in two sub-vectors of sizes $d_1, d_2$, such that $d_1+d_2=d$,
with corresponding partitions of the parameters in blocks of matching sizes,
as follows
\begin{equation}
 Y=\left( \begin{array}{c}
        Y_1 \\
        Y_2 \\
      \end{array}
    \right),\quad
 \xi=\left(  \begin{array}{c}
        \xi_1 \\
        \xi_2 \\
      \end{array}
    \right),\quad
 \gamma=\left( \begin{array}{c}
        \gamma_1 \\
        \gamma_2 \\
      \end{array}
    \right),
    \quad
 \Sigma  = \left(\begin{array}{cc}
                      \Sigma_{11} & \Sigma_{12} \\
                      \Sigma_{21} & \Sigma_{22} \\
                    \end{array}
        \right)\,.
  \label{e:partition}
\end{equation}
To establish the marginal distributions of $Y_1$, we use 
Proposition~\ref{th:affine} with $b=0$ and $B$ equal
to a matrix formed by $I_{d_1}$ in the top $d_1$ rows and
a block of $0$s in the bottom $d_2$ rows.
For $Y_2$, we proceed similarly, but setting the bottom $d_2$
rows of $B$ equal to $I_{d_2}$.
We then arrive at the following conclusion.

\begin{proposition} \label{th:marginal}
If $Y\sim\GMN_d(\xi, \Sigma, \gamma, H)$ is partitioned as indicated in
\ref{e:partition}, then
\begin{equation}
   Y_1\sim \GMN_{d_1}(\xi_1, \Sigma_{11}, \gamma_1, H), \qquad
   Y_2\sim \GMN_{d_2}(\xi_2, \Sigma_{22}, \gamma_2, H)\,.
   \label{e:gmn-marginals}
\end{equation}
\end{proposition}

We now want examine conditions which ensure independence of  $Y_1$ and $Y_2$.
From \ref{e:gmn-hierar} it is clear that, if $\Sigma_{12}=\Sigma_{21}\T=0$,
$Y_1$ and $Y_2$ are conditionally independent given $(U,V)$, 
with conditional distribution
\begin{equation}\label{e:gmn:marg-condi}
   (Y_j|U=u, V=v) \sim \N_d(\xi_j  + r(u,v)\,\gamma_j,\: s(u,v)^2\,\Sigma_{jj})\,,
   \quad j=1,2\,,
\end{equation}
where $(U,V) \sim G_U\times G_V.$ Moreover, if $s(U,V)\equiv 1$ (constant) and one of the
marginal distributions is symmetric,
i.e., $\gamma_1=0$ or $\gamma_2=0$, then $Y_1$ and $Y_2$  are independent.
The notation  $s(U,V)\equiv1$ and similar ones later on must be intended
`with probability 1'; we shall not replicate this specification subsequently.

A more detailed argument is as follows, where we take $\xi=0$ for mere simplicity of notation.
without affecting the generality of the argument.
From \ref{e:gmn-hierar},  we have that the conditional joint characteristic function 
of $(Y_1,Y_2)$, given $U=u$ and $V=v$ (or, equivalently, given $R=r$ and $S=s$), is
 \[
   \E{e^{it_1\T Y_1+it_2\T Y_2}\big|U=u,V=v}=
   \exp\left\{ir(u,v)(t_1\T\gamma_1+t_1\T\gamma_1)-
     \half s(u,v)^2(t_1\T\Sigma_{11}t_1+2t_1\T\Sigma_{12}t_2+t_2\T\Sigma_{22}t_2)\right\}
\]
so that the joint characteristic function of  $(Y_1,Y_2)$ is
 \begin{eqnarray}
 c(t_1,t_2)&=&\E{e^{it_1\T Y_1+it_2\T Y_2}}\nonumber\\
 &=&\E{\E{e^{it_1\T Y_1+it_2\T Y_2}|U,V}}\nonumber\\
 &=&\E{e^{ir(U,V)(t_1\T\gamma_1+t_2\T\gamma_2)-\half s(U,V)^2(t_1\T\Sigma_{11}t_1+2t_1\T\Sigma_{12}t_2+t_2\T\Sigma_{22}t_2)}}\nonumber\\
 &=&\E{e^{ir(U,V)t_1\T\gamma_1-\half s(U,V)^2t_1\T\Sigma_{11}t_1}e^{ir(U,V)t_2\T\gamma_2-\half s(U,V)^2t_2\T\Sigma_{22}t_2}e^{-s(U,V)^2t_1\T\Sigma_{12}t_2}}.
 \label{e:gmn:jcf}
 \end{eqnarray}
In analogous way, by \ref{e:gmn:marg-condi} the marginal characteristic functions are
\begin{eqnarray}
c_j(t_j)&=&\E{e^{it_j\T Y_j}}\nonumber\\
&=&\E{\E{e^{it_j\T Y_j}|U,V}}\nonumber\\
&=&\E{e^{ir(U,V)t_j\T\gamma_j-\half s(U,V)^2t_j\T\Sigma_{jj}t_j}},\quad j=1,2.\label{e:gmn:mcf}
\end{eqnarray}
Note that, if $\gamma_j=0$ and  $s(U,V)\equiv1$, then  by \ref{e:gmn:mcf} $c_j(t_j)$ reduces to the
centred normal characteristic function $c_{N,j}(t_j)=e^{-\half t_j\T\Sigma_{jj}t_j}$ for $j=1,2$.
We have then  reached the following conclusion.

\begin{proposition} \label{th:condition-independence}
Given partition \ref{e:partition}, the components $Y_1, Y_2$ are independent
provided $s(U,V)\equiv1$, $\Sigma_{12}=0$ and at least one of $\gamma_1$
and $\gamma_0$ is $0$, with the following implications:
\begin{itemize}
\item[~(a)] if $\gamma_1=0$, the joint characteristic function \ref{e:gmn:jcf}
   reduces to $c_{N,1}(t_1)\:c_2(t_2)$,
\item[~(b)]  if $\gamma_2=0$, the joint characteristic function \ref{e:gmn:jcf}
   reduces to $c_1(t_1)\:c_{N,2}(t_2)$.
\end{itemize}
If both $\gamma_1$ and $\gamma_2$ are $0$, the distribution reduces to the case
of independent normal variables.
\end{proposition}

In essence, under the conditions of Proposition~\ref{th:condition-independence},
one of $Y_1$ and $Y_2$  has a plain normal distribution and the other one falls
under the construction discussed later in Section~\ref{s:mean-mix}.

Outside the conditions of Proposition~\ref{th:condition-independence},
the structure of \ref{e:gmn:jcf} does not appear to be suitable for factorization
as the product of two legitimate characteristic functions, and we conjecture that,
in general,  independence between $Y_1$ and $Y_2$ cannot be achieved.

Examine now the conditional distributions associated to partition \ref{e:partition}.
Factorize the joint density of $Y$ as $f(y_1,y_2)=f_{1|2}(y_1|y_2)f_2(y_2)$
where $f_{1|2}(y_1|y_2)$ is the conditional density of $(Y_1|Y_2=y_2)$ and $f_2(y_2)$
is the marginal density of $Y_2$.
For simplicity of treatment, suppose that $(R,S)$ is absolutely continuous, with density $h(r,s)$.
Then, by \ref{e:pdf-gmn} and the properties of the multivariate normal density,  write
\begin{eqnarray*}
f_{1|2}(y_1|y_2)f_2(y_2) = \int_{\Real\times\Real_+}
  \phi_{d_1}(y_1;\xi_{1|2}+\gamma_{1|2}\,r, s^2\Sigma_{1|2})
  \phi_{d_2}(y_2;\xi_2+\gamma_2\,r, s^2\Sigma_{22})h(r,s) \d r\d s
\end{eqnarray*}
where
\[
 \xi_{1|2}=\xi_1+\Sigma_{12}\Sigma_{22}\inv(y_2-\xi_2), \quad
 \gamma_{1|2}=\gamma_1-\Sigma_{12}\Sigma_{22}\inv\gamma_2, \quad
 \Sigma_{11|2}=\Sigma_{11}-\Sigma_{12}\Sigma_{22}\inv\Sigma_{21}
\]
having assumed that the conditioning operation and integration can be exchanged.
Hence, for the conditional density of $Y_1$ given $y_2$ we have
\[
f_{1|2}(y_1|y_2) =  \frac{1}{f_2(y_2)}\:
  \int_{\Real\times\Real_+}\phi_{d_1}(y_1;\xi_{1|2}+\gamma_{1|2}\,r, s^2\Sigma_{1|2})
  \phi_{d_2}(y_2;\xi_2+\gamma_2\,r, s^2\Sigma_{22}) \:h(r,s)\:\d r\d s.
\]
Now, from the Bayes's rule, we obtain that the conditional density of $(R,S)$ given $Y_2=y_2$ is
\begin{equation}
  h_c(r,s|y_2)=\frac{\phi_{d_2}(y_2;\xi_2+\gamma_2\,r, s^2\Sigma_{22})\:h(r,s)}{f_2(y_2)}.
  \label{e:hc}
\end{equation}
Using this fact in the last integral, we can re-write
\begin{equation}
   f_{1|2}(y_1|y_2)= \int_{\Real\times\Real_+}
         \phi_{d_1}(y_1;\xi_{1|2}+\gamma_{1|2}\,r, s^2\Sigma_{11|2})\:h_c(r,s|y_2) \d{r}\d{s}
   \label{e:f1|2}
\end{equation}
which exhibits the same structure of \ref{e:pdf-gmn}.  Therefore we can conclude that
\[
  (Y_1|Y_2=y_2) \sim\GMN_{d_1}\left(\xi_{1|2}, \Sigma_{11|2}, \gamma_{1|2}, H_{c(y_2)}\right)
\]
where $H_{c(y_2)}$ denotes the distribution function associated to the conditional density \ref{e:hc}.

For many GMN constructions \ref{e:gmn}--\ref{e:gmn-RS}, the density function of $Y$
is likely  to be known in explicit form;
in these cases, the same holds true for $Y_2$, recalling \ref{e:gmn-marginals}.
Then, a convenient aspect of expression \ref{e:hc} is that it indicates how  to
compute the conditional density once the  joint unconditional distribution $H(r,s)$
is  available explicitly.
Clearly, this is especially amenable in those constructions where $(R, S)$ is really a
univariate variable, as in Sections~\ref{s:mean-mix} and~\ref{s:var-mix} below.

In one of the appendices,  we illustrate the use of \ref{e:hc}--\ref{e:f1|2}
in the case of a multivariate $t$ distribution.


\subsection{On quadratic forms}   \label{s:quadratic-forms}

For use in the next result, but also in the rest of the paper, define the quantities
\begin{equation}
   \Omega = \Sigma+\gamma\gamma\T, \quad
    \eta  = \left(1+\gamma\T\Sigma\inv\gamma\right)^{-1/2}\Sigma\inv\gamma\,,\quad
   \alpha^2= \|\gamma\|^2_\Sigma = \gamma\T\Sigma\inv\gamma,\,\quad
   \delta^2=\dfrac{\alpha^2}{1+\alpha^2}
   \label{e:aux-param}
\end{equation}
such that
$\alpha^2\in[0,\infty)$ and  $\delta^2\in[0,1)$.
For notational convenience, we introduce the notation
\begin{equation}
  \mu_{hk} = \E{R^h\; S^k},
  \qquad k=0,1,\dots   \label{e:mu.hk}
\end{equation}
when the named expectation exists.

\begin{proposition}
For a random variable $Y_0$ having distribution of type \ref{e:gmn}--\ref{e:gmn-RS}
with $\xi=0$, the following facts hold:
\begin{eqnarray}
  S^{-2}\, (Y_0-R\gamma)\T\Sigma\inv(Y_0-R\gamma)      &\sim& \chi_d^2,  \label{e:chi^2}\\
  \E{Y_0\T\Sigma\inv Y_0}  &=& d\:\E{S^2}+ \alpha^2\,\E{R^2}
        = d\:\mu_{02} +\alpha^2\mu_{20} \,,          \label{e:E(Mdist)}\\
  \E{Y_0\T\Omega \inv Y_0}
       &=&   d\:\E{S^2} + \delta^2\:\left(\E{R^2} - \E{S^2}\right)
          = d\:\mu_{02}+\delta^2(\mu_{20}-\mu_{02})   \,,   \label{e:E(Q)}
\end{eqnarray}
provided $\E{R^2}$ and $\E{S^2}$ exist, using the quantities defined in \ref{e:aux-param}
and \ref{e:mu.hk}.
\end{proposition}
\noindent{\it Proof:}
From \ref{e:gmn-RS}, write $(Y_0-R\gamma)\T\Sigma\inv(Y_0-R\gamma)=S^2\:X\T\Sigma\inv X$,
 where $X\T\Sigma\inv X\sim\chi_d^2$ is independent of $S$;
this yields result \ref{e:chi^2}.
For equality \ref{e:E(Mdist)}, expand the initial identity of this proof as
\[
 Y_0\T\Sigma\inv Y_0-2R\gamma\T\Sigma\inv Y_0+R^2\gamma\T\Sigma\inv\gamma=S^2\,X\T\Sigma\inv X
\]
and take expectation on both sides of this equality. We obtain
\[
  \E{R Y_0}=\E{R\:\E{Y_0|U,V}}
    =\E{R\:\E{R\,\gamma+S\,X)|U,V}}
    = \E{R\,(R\gamma+S\E{X|U,V})}=\E{R^2}\gamma,\]
bearing in mind that $\E{X|U,V}=\E{X}=0$, by the independence assumption
between $X$ and $(U,V)$. This leads to  \ref{e:E(Mdist)}.
\par
For \ref{e:E(Q)}, write $Q = Y_0\T\Omega\inv Y_0$, and
$\E{Q} = \tr\left(\Omega\inv\E{Y_0Y_0\T}\right)$.
Using \ref{e:E(Y)} and \ref{e:var(Y)}, we obtain
$\E{Y_0Y_0\T} = \var{Y_0} + \E{Y_0}\E{Y_0\T} = \E{R^2}\gamma\gamma\T+\E{S^2}\Sigma$,
so that
\[\E{Q} = \E{R^2}\gamma\T\Omega\inv\gamma+\E{S^2}\tr\left(\Omega\inv\Sigma\right).\]
By using the Sherman-Morrison equality \ref{e:SM}, we conclude the proof.
\hfill\textsc{qed}

\par\vspace{2ex}
In the subsequent pages, the matrix $\Omega$ defined in \ref{e:aux-param} and the
associated quadratic form $Q=Y_0\T\Omega\inv Y_0$ will appear repeatedly.
A connected relevant question is: under which conditions
is \ref{e:E(Q)}  free of $\gamma$? Equivalently, under which conditions
\begin{equation}
    \E{Q}= \E{Y_0\T\Omega \inv Y_0} =  d\:\E{S^2} \quad?   \label{e:E(Q)=dE(S)}
\end{equation}

This equality represents a form of invariance which is known to hold in
some cases to be recalled later on, but we want to examine  it
more generally.
One setting  where equality \ref{e:E(Q)=dE(S)} holds is given by
$R=U\,V^{1/2}$,  $S=V^{1/2}$, where $V>0$, and $\E{U^2}=1$.
It is then immediate to see that $\E{R^2}= \E{S^2}$,
so that the final term of \ref{e:E(Q)} is zero.

The conditions $R=U\,V^{1/2}$ and  $S=V^{1/2}$ are in turn achieved
when $Z=U\gamma+X$ and $Y_0=V^{1/2} Z$.
In this case $Q=VQ_0$, where $Q_0=Z\T\Omega\inv Z$ which is independent of $V$.
Hence, $\E{Q}=\E{V}\E{Q_0}$, where
\[\
    \E{Q_0}=\tr\left(\Omega\inv\E{Z\,Z\T}\right)
        =\E{U^2}\gamma\T\Omega\inv\gamma + \tr\left(\Omega\inv\Sigma\right)
\] since $\E{Z}=\E{U}\gamma$ and $\var{Z}=\var{U}\gamma\gamma\T+\Sigma$ and
so $\E{Z\,Z\T}=\E{U^2}\gamma\gamma\T+\Sigma$.
Thus, if $\E{U^2} = 1$, then by using \ref{e:SM} it clearly follows that $\E{Q_0}=d$.
We shall return to this issue later on.

\subsection{Mardia's measures of multivariate asymmetry and kurtosis}
\label{s:beta-Mardia}

For a multivariate random variable $Z$ such that $\mu_Z=\E{Z}$ and $\Sigma_Z=\var{Z}$,
\cite{mardia:1970,mardia:1974} has introduced measures of multivariate skewness and
kurtosis, defined as
\begin{equation}
  \beta_{1,d}=\E{\left[(Z-\mu_Z)\T\Sigma_Z\inv(Z'-\mu_Z)\right]^3} \,, \qquad
  \beta_{2,d}=\E{\left[(Z-\mu_Z)\T\Sigma_Z\inv(Z-\mu_Z)\right]^2},
  \label{e:beta-Mardia}
\end{equation}
where $Z'$ is an independent copy of $Z$, provided these expected values exist.
These measures represent extensions of corresponding familiar quantities 
for the univariate case:
\begin{equation}
  \beta_1 = \frac{\E{(Z- \mu_Z)^3}^2}{\var{Z}^3} = \gamma_1^2, \qquad
  \beta_2 = \frac{\E{(Z- \mu_Z)^2}}{\var{Z}^{2}} = \gamma_2 + 3,
  \label{e:beta-gamma-univar}
\end{equation}
in the sense that $\beta_{1,1}=\beta_{1}$ and $\beta_{2,1}=\beta_2$.

We want to find expressions for \ref{e:beta-Mardia} in the case of a 
random variable $Y$ of type \ref{e:gmn}--\ref{e:gmn-RS}.
Recall the expressions for $\mu_Y$ and $\Sigma_Y$ given in \ref{e:E(Y)} and \ref{e:var(Y)},
and the notation defined in \ref{e:mu.hk} for the moments of $(R,S)$, and write
\[ R_0=R-\mu_{10}\,, \qquad Y-\mu_Y=R_0\gamma+S\,X \,.\]
assuming that the involved mean values  exist.
Taking into account the invariance of $\beta_{d,1}$ and $\beta_{d,2}$
with respect to non-singular affine transformations, it is convenient to work
with the transformed quantities
\[
   X_0=\Sigma^{-1/2}\,X\sim \N_d(0,I_d), \qquad
  \gamma_0=\Sigma^{-1/2}\gamma,  \qquad
   Y_0=\Sigma^{-1/2}(Y-\mu_Y)=R_0\gamma_0+S\,X_0
\]
where any form of the square root matrix $\Sigma^{1/2}$ can be adopted.

The subsequent development involves extensive algebra of which we report here
only the summary elements;  detailed computations are provided in an appendix.
Recall $\alpha^2$ introduced in \ref{e:aux-param}  and  define 
\[
   \bar\mu_{20}=\mu_{20}-\mu_{10}^2=\var{R}\, \qquad
   \rho = \frac{\bar\mu_{20}}{\mu_{02}}\,, \qquad
   \bar\rho = \frac{\rho\alpha^2}{1+\rho\alpha^2}
            = \frac{\alpha^2 \bar\mu_{20}}{\mu_{02}+ \alpha^2 \bar\mu_{20}}\,.
\]
Introduce the auxiliary random variables $T_0=\alpha\inv\gamma_0\T X_0\sim \N(0,1)$,
which is independent of $(R,S)$, and $Z_0=\alpha\,R_0+S\,T_0$.
We need to compute the following expectations:
\begin{eqnarray*}
\E{S^2\,Z_0} &=&\alpha\:(\mu_{12}-\mu_{10}\mu_{02})\,,\\
\E{S^2\,Z_0^2}
  &=&\alpha^2(\mu_{22}-2\mu_{12}\mu_{10}+\mu_{10}^2\mu_{02})+ \mu_{04}\,,\\
\E{Z_0^3}
  &=&\alpha^3(\mu_{30}-3\mu_{20}\mu_{10}+2\,\mu_{10}^3) +
     3\alpha\:(\mu_{12}-\mu_{10}\mu_{02})\,,\\
\E{Z_0^4}
  &=&\alpha^4(\mu_{40}-4\mu_{30}\mu_{10}+6\mu_{20}\mu_{10}^2-3\,\mu_{10}^4)
     +6\alpha^2(\mu_{22}-2\mu_{12}\mu_{10}+\mu_{10}^2\mu_{02})+3\mu_{04}\,,
\end{eqnarray*}
assuming the existence of moments of $(R,S)$ up to the fourth order.
With these ingredients, the Mardia's measures for the GMN construction 
can  be expressed as
\begin{eqnarray}
\beta_{1,d} &=& \mu_{02}^{-3}
    \left(3(d-1)(1-\bar\rho)\:\E{S^2\,Z_0}^2+(1-\bar\rho)^3\:\E{Z_0^3}^2\right)\,,
  \label{e:beta1M-GMN}\\
\beta_{2,d} &=&\mu_{02}^{-2}
   \left((d+1)\,(d-1)\,\mu_{04}+2\,(d-1)\,(1-\bar\rho)\:\E{S^2\,Z_0^2}
        +(1-\bar\rho)^2\:\E{Z_0^4}\right)\,.
   \label{e:beta2M-GMN}
 \end{eqnarray}

Considering the complexity that typically involves the explicit specification
of \ref{e:beta-Mardia} outside the normal family, 
the above expressions appear practically manageable.
They are further simplified when one specializes them to a given family
or to a certain subclass of the GMN construction. For a given choice of
the distribution $H$, we need to work out the following ingredients:
(i)~the marginal moments of $R$, $\mu_{h0}$, up to order 4,
(ii)~the marginal moments $\mu_{02}$ and $\mu_{04}$ of $S$,
(ii)~the cross moments $\mu_{12}=\E{R\,S^2}$ and $\mu_{22}=\E{R^2\,S^2}$.
The working is illustrated next for the GH family; additional illustrations
will appear later.
\paragraph{Mardia's measures for the GH family} For the GH family with
representation \ref{e:var-mean-mix}, there is a single mixing variable
$V\sim \N^-(\lambda,\chi,\psi)$ with density \ref{e:pdf-gig} and
$R=V$, $S=V^{1/2}$.
General expressions for $\E{V^h}=\mu_{h0}$ are given in
Section~2.1 of \cite{jorgensen:1982}, among others. 
These expressions also provide $\mu_{02}=\E{V}$ and $\mu_{04}=\E{V^2}$.
The two other required quantities are $\mu_{12}=\E{V^2}$ and $\mu_{22}=\E{V^3}$
which are still ordinary moments of $V$.
We can now compute

\begin{eqnarray*}
\E{S^2\,Z_0}
&=&\alpha(\E{V^2}-(\E{V})^2)=\alpha\:\sigma^2_V\,, \hbox{~say},\\
\E{S^2 Z_0^2}
&=&\alpha^2(\E{V^3}-2\E{V^2}\E{V}+(\E{V})^2)+ \E{V^2},\\
\E{Z_0^3}
&=&\alpha^3(\E{V^3}-3\E{V^2}\E{V}+2(\E{V})^3) + 3\alpha\var{V}\\
&=&(\alpha\sigma_V)^3\beta_1(V) + 3\alpha\:\sigma^2_V \,,\\
\E{Z_0^4}
&=&\alpha^4(\E{V^4}-4\E{V^3}\E{V}+6\E{V^2}(\E{V})^2-3(\E{V})^4)\\
&&+6\alpha^2(\E{V^3}-2\E{V^2}\E{V}+(\E{V})^3)+3\E{V^2}\\
&=&(\alpha\sigma_V)^4\beta_2(V)+6\alpha^2(\E{V^3}-2\E{V^2}\E{V}+(\E{V})^3)+3\E{V^2}.
\end{eqnarray*}
where  $\sigma^2_V=\var{V}$ and $\beta_1(V)$, $\beta_2(V)$ are the univariate
measures of skewness and kurtosis in \ref{e:beta-gamma-univar} evaluated for $V$.
Plugging the above quantities in \ref{e:beta1M-GMN} and \ref{e:beta2M-GMN}
completes the computation.

\paragraph{Remark} 
There exists an interesting way of re-writing \ref{e:beta1M-GMN} and \ref{e:beta1M-GMN} 
which will turn out useful later on.
 Since $Z_0=\alpha\,R_0+S\,T_0\sim \GMN_1(-\alpha\mu_{10},1,\alpha,H))$ 
with  zero mean and  
\[\var{Z_0}=\alpha^2\bar\mu_{20}+\mu_{02}=(1-\bar\rho)\inv\mu_{02}\] 
we can introduced an univariate standardized GMN-type variable 
\[  
\tilde Z_0=\frac{\alpha\,R+S\,T_0-\alpha\,\mu_{10}}{\sqrt{\alpha^2\bar\mu_{20}+\mu_{02}}}\sim 
\GMN_1\left(-\frac{\mu_{10}}{\sqrt{\alpha^2\bar\mu_{20}+\mu_{02}}},
         \frac{1}{\alpha^2\bar\mu_{20}+\mu_{02}},
          \frac{\alpha}{\sqrt{\alpha^2\bar\mu_{20}+\mu_{02}}},H\right)
\] 
which has zero mean zero and unit variance. 
When rewritten in terms of $\tilde Z_0$,  
\ref{e:beta1M-GMN} and \ref{e:beta2M-GMN} become
\begin{eqnarray}
\beta_{1,d} &=& 3(d-1)\mu_{02}^{-2}\:\E{S^2\,\tilde{Z}_0}^2+ \beta_1(\tilde Z_0) \,,
      \label{e:beta1M-GMN-Z0}  \\
\beta_{2,d} &=& (d+1)\,(d-1)\,\mu_{02}^{-2}\mu_{04}+2\,(d-1)\,\mu_{02}\inv\:\E{S^2\,\tilde{Z}_0^2}
        + \beta_2(\tilde Z_0)\,,
       \label{e:beta2M-GMN-Z0}   
\end{eqnarray}
where $\beta_1(\tilde Z_0)$ and $\beta_2(\tilde Z_0)$ denote the univariate coefficients 
$\beta_1$ and $\beta_2$ in \ref{e:beta-gamma-univar} evaluated for $\tilde Z_0$.

\section{Mean (or location) mixtures}  \label{s:mean-mix}
In this section and the next one, we examine two simplified versions of the
general formulation  \ref{e:gmn}--\ref{e:gmn-RS}.
The first class occurs when only the additive random component is actually
present.
\subsection{General properties}   \label{s:mean-mix-general}

\paragraph{Basic facts}
Consider the simplified form of \ref{e:gmn}--\ref{e:gmn-RS} where $s(u,v)=1$,
so that   we can assimilate $R$ and $U$, and write
\begin{equation}
  Y = \xi + U\,\gamma + X \,.
  \label{e:mean-mix}
\end{equation}
In this case we use the notation $Y\sim\GMN_d(\xi, \Sigma, \gamma, G_U)$,
since now $(R,S)$ reduces to $R\equiv U$, with distribution $H=G_U$.

Clearly, if $U$ is a degenerate random variable, $U\equiv1 $ say,  or 
if $\gamma=0$, the construction reduces to the normal distribution.
We therefore exclude these cases from consideration.

Although constructions of this type can simply be viewed as a sum of two
independent random components, they can legitimately also be
regarded as a location mixture, within the logic of \ref{e:gmn-hierar},
and this interpretation  can facilitate the construction of EM-type
algorithms and work in Bayesian inference.

Several general properties of the class \ref{e:mean-mix} have been obtained
by \cite{nega:jama:etal:2019}.
Their initial development is in the univariate context, where they establish
the property of closure  under convolution with an independent normal variables,
and an expression of the characteristic function. From these results, they
derive expressions for low order moments and associated measures of
asymmetry and kurtosis.
Section~8 of their paper refers to the multivariate case, where they
also obtain various results, notably  the property of closure under
marginalization, an expression for the conditional distribution given
the values taken on by certain components of $Y$, and
\begin{equation}
   \E{Y} = \xi + \E{U}\,\gamma\,, \qquad \quad
   \var{Y} =  \var{U}\:\gamma\gamma\T + \Sigma
\label{e:gmn-loc,E&var}
\end{equation}
provided $\E{U}$ and $\var{U}$ exist. These expressions can also
be obtained as special cases of \ref{e:E(Y)} and \ref{e:var(Y)}.

The property of closure under convolution with normal variates,
which has been stated by  \cite{nega:jama:etal:2019} in the univariate
case, actually holds also in the multivariate case.
Specifically, if $Y\sim\GMN_d(\xi, \Sigma, \gamma, G_U)$  and
$W\sim\N_d(\mu,\tilde\Sigma)$ are independent variables, then
it is immediate from representation \ref{e:mean-mix} that
$Y+W\sim\GMN_d(\xi+\mu, \Sigma+\tilde\Sigma, \gamma,  G_U)$.

From Proposition~\ref{th:condition-independence},
we can say that the marginal components $Y_1$  and $Y_2$ are independent if
and only if $\Sigma_{12}=0$ and at least one of
$\gamma_1$ and $\gamma_2$ is $0$.


\paragraph{Mardia's measure for mean mixtures}
With respect to development in Subsection~\ref{s:beta-Mardia}, here we have
$S\equiv1$ and $R=U$, with substantial simplification of the general
expressions \ref{e:beta1M-GMN}--\ref{e:beta2M-GMN}. 
In this case we obtain
\begin{eqnarray*}
\E{S^2 Z_0} &=& \E{Z_0} = 0,\\
\E{S^2 Z_0^2} &=& \E{Z_0^2}= (\alpha\sigma_U)^2 + 1=(1-\bar\rho)\inv\\
\E{Z_0^3} &=& \alpha^3(\mu_{30}-3\mu_{20}\mu_{10}+2(\mu_{10})^3)\\ 
  &=& (\alpha\sigma_U)^3\,\beta_1(U),\\
\E{Z_0^4} &=& \alpha^4(\mu_{40}-4\mu_{30}\mu_{10}+6\mu_{20}(\mu_{10})^2-3(\mu_{10})^4)+6\alpha^2\sigma_U^2+3\\
&=&(\alpha\sigma_U)^4\beta_2(U)+6(\alpha\sigma_U)^2+3.
\end{eqnarray*}
where $\beta_1(U)$ and $\beta_2(U)$ denote the univariate coefficients in
\ref{e:beta-gamma-univar} evaluated for $U$, leading to
\begin{eqnarray*}
\beta_{1,d} &=& (1-\bar\rho)^3(\E{Z_0^3})^2\nonumber\\
&=&((\alpha\sigma_U)^2+1)^{-3}[(\sigma_U\alpha)^3\,\beta_1(U)]^2,
 \\
\beta_{2,d} &=&d(d+2)-3+(1-\bar\rho)^2\E{Z_0^4}\nonumber\\
&=&d(d+2)-3+((\alpha\sigma_U)^2+1)^{-2}[(\sigma_U\alpha)^4\,\beta_2(U)+6(\alpha\sigma_U)^2+3]\nonumber\\
&=&d(d+2)+((\alpha\sigma_U)^2+1)^{-2}[-3((\alpha\sigma_U)^2+1)+(\sigma_U\alpha)^4\,\beta_2(U)+6(\alpha\sigma_U)^2+3]\nonumber\\
&=&d(d+2)+((\alpha\sigma_U)^2+1)^{-2}[(\alpha\sigma_U)^4\beta_2(U)+3(\alpha\sigma_U)^2] \,.
\end{eqnarray*}

Note that the leading term $d(d+2)$ in the last expression represents the difference
between $\beta_{2,d}$ and its companion measure of excess, $\gamma_{2,d}$,
in \cite{mardia:1974}.

\paragraph{Remark} In this case, \ref{e:beta1M-GMN-Z0} and \ref{e:beta2M-GMN-Z0} 
yield a very neat simplification, namely
\begin{eqnarray*}
\beta_{d,1} &=& \E{\tilde{Z}_0^3}^2=\gamma_1^2\nonumber\\
\beta_{2,d} &=&d(d+2)+\E{\tilde{Z}_0^4}-3=d(d+2)+\gamma_2\nonumber\,.
\end{eqnarray*}

\subsection{Some noteworthy special cases}   \label{s:mean-mix-classes}

The more interesting families of this class are arguably those obtained
when  the distribution of $U$ is not symmetric about 0.
In fact, in nearly all special formulations discussed below,
$U$ is a positive variable.

Besides its intrinsic values from the distribution theory viewpoint,
there is the interesting connection of \ref{e:mean-mix}  with non-symmetric
$U$ and the formulation
in quantitative finance put forward  by \cite{simaan:1993}, as for the
assumptions on the key stochastic component, and the closure under
marginalization.

\paragraph{The skew-normal distribution and its extended version}
When $U$ in \ref{e:mean-mix} has a positive half-normal distribution,
or equivalently the $\N(0, 1)$ distribution truncated below $0$,
we obtain the set-up adopted by \citet{azza:dval:1996} to derive
the density function of the multivariate skew-normal (SN) family.
The multivariate SN density function  at $y\in\Real^d$  is
\begin{equation}
  2\: \phi_d\left(y-\xi; \Omega\right)\:\Phi\left(\eta\T(y-\xi)\right)
\label{e:pdf-sn}
\end{equation}
where $\Omega$ and $\eta$  are as in \ref{e:aux-param}.
In one appendix, we present a  proof of this expression which
retains the same logic of the proof of \citet{azza:dval:1996},
but involves a more essential development.

The multivariate SN  distribution enjoys a number of appealing
formal properties, matching many of those of the normal distribution.
An account of this theme is provided in Chapter~5 of \citet{azza:capi:2014}.
In view of the discussion in  Section~\ref{s:quadratic-forms},
we must at least mention  the fact that, as a special case of
a more general result on quadratic forms,
$(Y-\xi)\T\Omega\inv(Y-\xi)\sim\chi^2_d$ when $Y$ is a random
variable with density \ref{e:pdf-sn}.

The extended form of the skew-normal distribution occurs when $U$
is distributed as $\N(0, 1)$ variable truncated below $-\tau$ instead
of $0$, for some constant $\tau$. The distribution of $U$ is then
$\phi(u)/\Phi(u)$ for $u+\tau>0$. A simple adaptation of the
above-mentioned proof yields the density function of $Y$
at $y\in\Real^d$  as
\begin{equation}
  \frac{1}{\Phi(\tau)} \:
  \phi_d\left(y-\xi; \Omega\right)\:\Phi\left(\bar\tau + \eta\T(y-\xi)\right)
\label{e:pdf-esn}
\end{equation}
where $\bar\tau=\left(1+\gamma\T\Sigma\inv\gamma\right)^{1/2}\:\tau$. 

While it is not clear whether a EM-type approach is the most efficient
way to tackle maximum likelihood estimation for distribution
\ref{e:pdf-sn} or \ref{e:pdf-esn}, certainly
EM-type  algorithms based on \ref{e:gmn-hierar} constitute
a popular route for parameter estimation in this context.
Early publications adopting this route to estimation include
\cite{arel:bolf:lach:2005} and \cite{arel:ozan:etal:2005},
but many others exist, often in connection with finite mixtures
of SN distributions.

\paragraph{MMNE and MMMNE distributions}\,
A substantial portion of the paper of \cite{nega:jama:etal:2019}
focuses on the specific instance where $U$ in \ref{e:mean-mix}
follows a standard exponential distribution.
They initially examine the case where $Y$ is univariate; this is said to have
a MMNE distribution, and  several interesting properties are derived:
log-concavity of the density, monotonicity of the hazard rate, infinitely
divisibility and more.
Subsequently, they consider the multivariate version,
called MMMNE distribution, whose density function at $y\in\Real^d$
can be written, with an inessential notational variation from the original paper,
as
\begin{equation}
  \phi_d(y;\xi,\Sigma)\:\left\{\left(\gamma\T\Sigma\inv\gamma\right)^{1/2}
      \hspace{1ex}
    \zeta\left(\frac{\gamma\T\Sigma\inv{x}-1}%
        {~\left(\gamma\T\Sigma\inv\gamma\right)^{1/2}}\right)
    \right\}\inv
\end{equation}
where $\zeta(\cdot)$ is defined in \ref{e:zeta}.

Additional results are derived for the MMMNE distribution, notably
the expression of the characteristic functions, the marginal and
the conditional distributions given the value taken on by a subset of $Y$
components. The mean and the variance are simply obtained by setting
$\E{U}=\var{U}=1$ in  \ref{e:gmn-loc,E&var}.

\paragraph{A two-piece normal mixing} Assume that $U$ in \ref{e:mean-mix} has a two-piece
normal distribution, that is, one having density function at $u\in\Real$:
\[
    2\pi_a\phi(u;a^2)I_{(0,\infty)}(u)+2\pi_b\phi(u;b^2)I_{(-\infty,0]}(u)
\]
where $I_A(\cdot)$ denotes the indicator function of set $A$.
This construction has been proposed repeatedly in the literature as a
simple way to allow for skewness via a simple modification of the normal density.
A compilation of rediscoveries of this distribution has been presented by
\cite{wallis:2014}.

It has been shown by \cite{arel:azza:etal:2020} that $Y$, as defined in \ref{e:mean-mix},
has a density function represented by a two-component mixture of skew-normal
variates. Hence each component has a density of type \ref{e:pdf-esn} with
$\tau=0$.

This distribution of $U$ is the only instance reviewed here where
$U$ is not a positive variable. It is included in our list
because it represents an interesting bridge between different families:
it shows how a mixture of multivariate normal variates, suitably combined
with a mixing two-piece normal variate, yields a multivariate SN variable.
In other words, it provides a link between two asymmetric extensions of the
normal family, the two-piece and the skew-normal distributions.


\paragraph{A Rayleigh mixing mean}
As far as we know, the following construction has not been examined
in the literature.
Suppose that the mixing variable $U$ in \ref{e:mean-mix} has a standard Rayleigh
distribution, with density
\begin{equation}
   g_U(u) = u e^{-\half u^2}I_{[0,\infty)}(u).
   \label{e:pdf-rayleigh}
\end{equation}
and write $Z=U \gamma + X$.
In this case, the mean and the variance of $Y=\xi+Z$ are simply obtained by setting
$\E{U}=(\pi/2)^{1/2}$ and $\var{U}=(4-\pi)/2$ in  \ref{e:gmn-loc,E&var}.

Consideration of this model can be motivated as follows. 
The Rayleigh distribution is widely used in a various applied disciplines,
especially in the engineering context. An example is represented by the 
technology of wind energy, where the distribution  plays a central role, 
as clearly visibile in the comprehensive treatment
of the subject by  \citet{nelsonV:2013}. 
From this source, we underline the noteworthy fact that ``Manufacturers 
[of wind turbines] assume a Rayleigh distribution for a wind speed'' (p.\,101).
Besides this domain,  distribution \ref{e:pdf-rayleigh} is used in
signal processing, ocean energy  and off-shore engineering, 
material design and reliability, and other areas, not all in engineering.
Consider now the case where a certain event, such as wind speed $U$ at a
certain location and time, is measured by $d$ instruments at the time, 
not just one. There will then be a $d$-dimensional vector $X$ 
of random components originated by the measuring instruments
which, once combined with $U$,  yields an instance of \ref{e:mean-mix}.
If all the instruments are perfectly calibrated, $\xi$ will be
the null vector and $\gamma$ will have all components equal,
otherwise discrepancies will exists.

It will be noted that density \ref{e:pdf-rayleigh} does not include 
a positive scale factor, $\sigma$ say, which is essential in applied work.
This factor is implicit and subsumed in $\gamma$, as otherwise
we would incur in a overparameterization situation.
In an applied context, it would presumably be sensible to
reparameterize in some more meaningful form, 
such as $\gamma=\sigma\tilde\gamma$,
with some suitable constraint on $\tilde\gamma$.
In the present more technical context, we retain the use of $\gamma$.

We show in one appendix that the density function of $Z$ is
\begin{equation}
  f_Z(z; \gamma, \Sigma) =
      (2\pi)^{1/2}(1+\alpha^2)^{-1/2}\phi_d(z;\Omega)\:\Phi(\eta\T z)\:
     \left\{\zeta(\eta\T z)  + \eta\T z\right\}
  \label{e:pdf-Rayleigh-mix}
\end{equation}
where $\zeta(\cdot)$ is defined in \ref{e:zeta}.
From \ref{e:pdf-Rayleigh-mix}, it is immediate to obtain the density of $Y$.

It is interesting that the leading factor on the right-hand side of \ref{e:pdf-Rayleigh-mix}
is, up to a constant, the SN density \ref{e:pdf-sn} with $\xi=0$.
In the specific case with $d=1$ and $\Sigma=1$, this factor is the $\SN(0,1+\gamma^2,\gamma)$ density.
 
\paragraph{An extension to $\chi_\nu$ mixing} \,
In two of the constructions examined earlier,
the $U$ component is a square-root of a $\chi^2$ variable.
Specifically, in the SN construction, $U$ has as a half-normal distribution,
that is, $U\sim\chi_1$. In another case examined, $U$ has a Rayleigh
distribution, that is, $U\sim\chi_2$.
It is then natural to consider a more general formulation where
$U\sim \chi_\nu$ for some positive $\nu$, having density
\[
  g_U(u)=\frac{2(1/2)^{\nu/2}}{\Gamma(\nu/2)} \: u^{\nu-1} e^{-\half u^2}, \qquad
  u\in\Real^+\,.
\]

The remark made in connection with \ref{e:pdf-rayleigh} about incorporation of
any scale parameter of $U$ in $\gamma$ carry on here.
In the light of this, the $\chi_\nu$ distribution considered here is effectively 
equivalent to what is often called Nakagami $m$-distribution \citep{nakagami:1960} 
in radio communication engineering.

Recalling the expression of moments of the $\chi^2$ distribution, we obtain readily
\[
   \E{U^m}=\frac{2^{m/2}\Gamma((\nu+m)/2)}{\Gamma(\nu/2)}, \qquad \var{U}=2\nu-\E{U}^2.
\]
The density of $Z=\gamma U+ X$ now becomes
\begin{eqnarray}
f_Z(z;\nu,\gamma,\Sigma)
  &=& \frac{2\sqrt{\pi}}{\Gamma(\nu/2)[2(1+\alpha^2)]^{(\nu-1)/2}} \:
     \phi_d(z;\Omega)\int_{-\eta\T z}^\infty(w + \eta\T z)^{\nu-1}\phi(w)\d w 
      \nonumber \\
  &=& \frac{2\sqrt{\pi}}{\Gamma(\nu/2)[2(1+\alpha^2)]^{(\nu-1)/2}} \: \phi_d(z;\Omega)
    \:\Phi(\eta\T z)  \: M_{\nu-1}(\eta\T\zeta)
    \label{e:pdf-chi-mix}
\end{eqnarray}
having written
\[
  M_{\nu-1}(\eta\T z)=  \E{(W + \eta\T z)^{\nu-1}\mid W + \eta\T z>0}
\]
where $W\sim \N(0,1)$.

To compute $ M_{\nu-1}(\eta\T z)$, we  must assume that $\nu$ is integer.
We then expand
\[
  M_{\nu-1}(\eta\T z)
  = \sum_{k = 0}^{\nu-1}
   {{\nu -1} \choose k} \:\E{W^{\nu-1}\mid W + \eta\T z>0}(\eta\T z)^{\nu-k-1}
\]
where $\E{W^{\nu-1}\mid W + \eta\T z>0}$ is the $(\nu-1)$th moment of a truncated
normal distribution.

A recursive expression for the moments of a $\N(0,1)$ variable truncated below
level $a$, say, has been given by \cite{elandt:1961}, which in our case must
be used with $a=-\eta\T z$.
An alternative route to these moments is via consideration of the
moment generating function and the cumulant generating
function of the truncated normal distribution, namely,
\[
   M(t)=e^{t^2/2}\,\frac{\Phi(t-a)}{\Phi(-a)},\qquad
   K(t)=\half t^2+\log\{\Phi(t-a)\}-\log\{\Phi(-a)\}
\]
whose derivatives involve those of $\zeta(\cdot)$ in \ref{e:zeta}. 
Expressions of low order derivatives of $\zeta(\cdot)$ are given in Section~2.1.4 
of \cite{azza:capi:2014}.

To develop a EM-type algorithm for this distribution, start by considering
the following hierarchical representation:
\begin{eqnarray*}
  Y | U &\sim& \N_d(\xi + \gamma U, \Sigma),\\
  U &\sim & \chi_\nu.
\end{eqnarray*}
Next, we must compute the first two moments of $\E{U|Y}$,
which in turn requires the conditional density of $U$ given $Y$ or, equivalently, given $Z$.
Taking into account that
\[\phi_d(z;\gamma u,\Sigma)=\phi_d(y;0,\Sigma)e^{-\half(\alpha^2 u^2-2\gamma\T\Sigma\inv z)}
=\phi_d(z;0,\Omega)e^{-\half\{\alpha^2 u^2-2\gamma\T\Sigma\inv z+(1+\alpha^2)\inv(\gamma\T\Sigma\inv z)^2\}},
\]
a standard application of Bayes theorem gives the conditional density
\begin{eqnarray*}
h_c(u|z)
&=&\frac{\phi_d(z;\gamma u,\Sigma)h(u)}{f_Z(z;\nu,\gamma,\Sigma)}\\
&=&\frac{\phi_d(z;0,\Omega) \:\exp\left(-\half\{\alpha^2 u^2-2\gamma\T\Sigma\inv z+(1+\alpha^2)\inv(\gamma\T\Sigma\inv z)^2\}\right)
\dfrac{2(1/2)^{\nu/2}}{\Gamma(\nu/2)} \: u^{\nu-1} e^{-\half u^2}}%
{\dfrac{2\sqrt{\pi}}{\Gamma(\nu/2)[2(1+\alpha^2)]^{(\nu-1)/2}} \: \phi_d(z;\Omega)\Phi(\eta\T z)
    \:M_{\nu-1}(\eta\T z)}\\
     &=&\frac{(1+\alpha^2)^{(\nu-1)/2} u^{\nu-1}
        \exp\left(-\half(1+\alpha^2)\{u-(1+\alpha^2)^{-1/2}\eta\T z\}^2\right)}%
        {\sqrt{2\pi} \:\Phi(\eta\T z) \: M_{\nu-1}(\eta\T z)}\\
    &=& \frac{(1+\alpha^2)^{(\nu-1)/2}}{M_{\nu-1}(\eta\T z)\: \Phi(\eta\T z)}\:
     \: u^{\nu-1}\phi\left[(1+\alpha^2)^{1/2}\left(u-(1+\alpha^2)^{-1/2}\eta\T z\right)\right],\quad u>0.
     %
\end{eqnarray*}
Computation of the $k$th moment of $h_c(u|z)$ effectively amounts to compute the
$(k+\nu-1)$th moment of a truncated normal distribution, 
a point which we have discussed earlier.

In this construction, we have left unspecified whether $\nu$ represents 
a known constant or a free positive parameter to be estimated. 
Similarly to the cases with $\nu=1$ and $\nu=2$ which correspond to
already-examined distributions, \ref{e:pdf-chi-mix}  could be employed 
with a fixed value of $\nu$, a situation which would ease use of the 
EM algorithm introduced above. 
There is, however, no bar to use it also when $\nu$ is a free integer parameter. 

\subsection{Again about quadratic forms}  \label{s:quadratic-forms2}

Consider the question of equality \ref{e:E(Q)=dE(S)},
and more generally the distribution of the quadratic form $Q$,
in the framework of mean mixtures \ref{e:mean-mix}.
Since now $R=U$ and $S \equiv 1$,  we consider the distribution of
$Q_0=Z\T\Omega\inv Z$, where $Z=U\gamma+X$.
A point of special interest are the conditions on $U$ such that $Q_0$
is distributed as $X\T\Sigma\inv X\sim\chi^2_d$, and $Q_0$ is independent of $(U,V)$.

Note that, if the distribution of $Q_0$ does not depend on $\gamma$,
the same holds true for the more general case where  $R=V^{1/2}\,U$ and $S=V^{1/2}$,
such that $Q=S^2\,Q_0$.
This setting falls within the more general construction examined in
Section~\ref{s:var-mean-mix}.

In the following discussion, we can ignore the special case $\gamma=0$,
as otherwise we return to the basic setting  of a normal variable $Z$
for which it is well-known that $Q_0 \sim \chi^2_d$.
Define $W_0=\gamma\T\Sigma\inv X\sim\N(0,\alpha^2)$, using the notation in \ref{e:aux-param};
note that $\alpha^2>0$.
On setting $X_0 = \Sigma^{-1/2}X$, we can also write  $W_0=\gamma\T\Sigma^{-1/2} X_0$.
For definitiveness, we take $\Sigma^{1/2}$ to be unique symmetric positive-definite
square root matrix of $\Sigma$, although we next steps would hold also with other
choices of the square root.
For notational convenience,  introduce $T_0 = \alpha\inv\,W_0 \sim\N(0,1)$.

\begin{proposition}  \label{th:Q0=W.sq+V0.sq}
Under the above definitions of symbols, $Q_0$ can be decomposed as
\begin{equation}
   Q_0 = W^2  + V_0^2
   \label{e:Q0=W.sq+V0.sq}
\end{equation}
where $W = \delta U  + (1-\delta^2)^{1/2} T_0$ and $V_0^2=\|X_0\|^2 - T_0^2$
are independent variables, with $V_0^2\sim\chi^2_{d-1}$.
\end{proposition}

\noindent{\it Proof:}
By applying the Sherman-Morrison formula \ref{e:SM} to $\Omega\inv$, we can write
\[
  \Omega\inv\gamma = (1 + \gamma\T\Sigma\inv\gamma)\inv\Sigma\inv\gamma,  \qquad
 \gamma\T\Omega\inv\gamma=(1 + \gamma\T\Sigma\inv\gamma)\inv\gamma\T\Sigma\inv\gamma
  = \delta^2
\]
so that we can decompose $Q_0$ as
\begin{eqnarray*}
Q_0 &=& (U\gamma + X)\T\Omega\inv(U\gamma + X)\\
 &=& U^2 \gamma\T\Omega\inv\gamma + 2 U\gamma\T\Omega\inv X + X\T\Omega\inv X\\
&=& U^2 (1 + \gamma\T\Sigma\inv\gamma)\inv\gamma\T\Sigma\inv\gamma + 2 U(1 + \gamma\T\Sigma\inv\gamma)\inv\gamma\T\Sigma\inv X
+ X\T\Sigma\inv X-(1 + \gamma\T\Sigma\inv\gamma)\inv(\gamma\T\Sigma\inv X)^2\\
&=& \delta^2 U^2  + 2\delta (1-\delta^2)^{1/2}U T_0 - \delta^2 T_0^2 + \|X_0\|^2\\
&=& \delta^2 U^2  + 2\delta (1-\delta^2)^{1/2}U T_0 + (1 - \delta^2) T_0^2 - T_0^2 + \|X_0\|^2\\
&=& \left(\delta U  + (1-\delta^2)^{1/2} T_0\right)^2  - T_0^2 + \|X_0\|^2
\end{eqnarray*}
which proves equality \ref{e:Q0=W.sq+V0.sq}.

\par

On defining the unit-norm vector $\bar\gamma = \alpha\inv\Sigma^{-1/2}\gamma\in\Real^d$,
we can write $V_0^2=X_0\T M_0 X_0= \|M_0 X_0\|^2$  where $M_0=I_d-\bar\gamma\bar\gamma\T$
is a symmetric idempotent matrix of rank $d-1$.
Hence, by standard results in normal theory distribution,
this proves the claim  that $V_0^2 \sim\chi^2_{d-1}$.
Moreover,  $T_0=\bar\gamma\T X_0$ and $M_0 X_0$ are orthogonal projections of $X_0$,
since $\bar\gamma\T M_0=0$, and then independent normal variables.
This implies independence of $T_0$ and $V_0^2$ and, since $X_0$ and its transformations
such as $T_0$ are independent of $U$, we conclude that $W$  and $V_0^2$ are independent.
\hfill\textsc{qed}
\par\vspace{2ex}

A corollary of Proposition~\ref{th:Q0=W.sq+V0.sq}, taking into account the  independence
of $U$ and  $T_0$, is that
\begin{eqnarray*}
  \E{Q_0} &=& \delta^2\:\E{U^2} + (1-\delta^2) + d-1  \\
          &=& d + \delta^2\left(\E{U^2}-1\right)
\end{eqnarray*}
provided $\E{U^2}$ exists.
Therefore, $\E{Q_0}$ does not depend on $\delta^2$, hence on $\gamma$,
if and only if $\E{U^2}=1$, in which case $\E{Q_0}=d$.

When $U$ is distributed as a positive half-normal, $W$ in \ref{e:Q0=W.sq+V0.sq} has
a univariate SN distribution.
By a known property of the SN distribution, we can say that $W^2\sim\chi^2_1$ and,
by using Proposition~\ref{th:Q0=W.sq+V0.sq}, we conclude that $Q_0\sim\chi^2_d$.
This is a well-known fact in the pertaining literature, as we have recalled  after
introducing  density  \ref{e:pdf-sn}. The derivation here is unusual, because
it is aimed to address the following question:
are there other choices of $U$ leading to the same distribution of $Q_0$?%

\section{Variance (or scale) mixtures}   \label{s:var-mix}

\subsection{General points}

In a sense, the dual formulation of mean mixtures is represented by variance
(or scale) mixtures, already mentioned in the introductory section.
This class can be produced by setting $r(u,v)=1$ in \ref{e:gmn} and subsuming
$\gamma$ into $\xi$, or equivalently by setting $\gamma=0$.
In either case, we arrive at formulation
\begin{equation}
  Y = \xi + V^{1/2}\: X
  \label{e:var-mix-2},
\end{equation}
 where we
have assimilated $S$ and $V$, with the condition $V>0$. The notation introduced for the
general construction \ref{e:gmn}--\ref{e:gmn-RS} now simplifies to
$Y\sim\GMN(\xi, \Sigma, 0, G_V)$, since $(R,S)$ reduces to $S=V^{1/2}$ 
 and $H$ reduces to $G_V$.
  
As already recalled in Section~\ref{s:normal-mix}, a scale  mixture of
normal distributions is a member of the class of elliptically contoured distributions.
In many popular constructions, $V$ is a continuous variable,
and this is the situation on which we shall focus.

In a large number of cases, the inverse operation is also possible,
that is, many members of the elliptical class can be represented as
scale mixtures of normal distributions. Note, however, that the implied
mixing variable $V$ has often a distribution which depends on the
dimension $d$ of $X$. This situation prevents the property of closure under
marginalization; on this issue, see \cite{kano:1994}.

We now recall some instances of scale mixtures of normal variables,
but only very briefly and confining ourselves to a few key instances,
since they represent very familiar constructions,
discussed in many existing accounts.
A classical early reference on this theme is \cite{andr:mall:1974}.
Several instances of the construction are presented by \citet{lang:sins:1993};
note that they use the alternative term `normal/independent distributions'
to identify this theme.
A relatively more recent account  is provided in Section~3.2
of \citet{mnei:frey:embr:2005}.

The formulation \ref{e:var-mix} bears the danger of over-parameterization, since
one could manoeuvre scale both via a suitable parameter of $V$ and via
the $\Sigma$ matrix. The issue is usually solved by ruling out  a scale parameter
in the distribution of $V$.

\subsection{Some noteworthy special cases}

\paragraph{Student's $t$ distribution}
Presumably, the best known parametric family of this class is
the Student's $t$ distribution, which occurs when $V\sim \nu/\chi^2_\nu$,
in an obvious notation.
In this case we write $Y\sim t_d(\xi, \Sigma,\nu)$.
Note that the distribution of $V$ does not allow for a scale parameter,
hence avoiding the above-mentioned issue of lack of identifiability.
For later use, recall the expression of the multivariate $t$ density function:
\begin{equation}
  t_d(y; \xi, \Sigma, \nu) =  \frac{\Gamma((\nu+d)/2)}%
        {(\nu\pi)^{d/2} \,\Gamma(\nu/2)\,\det(\Sigma)^{1/2}}
        \left(1+\nu\inv\:\|y-\xi\|^2_\Sigma\right)^{-\frac{\nu+d}{2}}  \,,
   \hspace{3em} y\in\Real^d,
  \label{e:pdf-t}
\end{equation}
where the expression $\|y-\xi\|^2_\Sigma$ makes use of the notation \ref{e:M-dist}.

An important property of this family is that all marginal distributions of $Y$
are still of Student's $t$ type  with the same degrees of freedom and
the other parameters as indicated in the partition \ref{e:partition}.
 This result hinges of the fact that
the distribution of $V$ does not depend on the dimension $d$.
The conditional distribution presented in one appendix shows that a
similar property holds also for the conditional distribution,
although in this case the degrees of freedom and other parameters
are  modified.
A wealth of other results on the multivariate $t$ distribution and its
variants or extensions is presented in the monograph of \citet{kotz:nada:2004}.

\paragraph{Symmetric GH distribution}
Consider the case where $V$ has a generalized inverse Gaussian
distribution \ref{e:pdf-gig}.
The implied GMN density of $Y$ is  the symmetric GH density
whose expression can be obtained by the general GH density shown
below in equation \ref{e:pdf-gh} when $\gamma=0$.

\paragraph{Other instances}
Of the many other instances of construction \ref{e:var-mix-2}, 
a mention is due for the case where $V$ is  discrete.
The basic instance is a two-point distribution, leading to
the contaminated normal distribution.

Another interesting instance is the multivariate slash distribution which
arises choosing $V=W^{-r}$ where $W\sim U(0,1)$ and $r$ is a positive value,
which regulates the tail weight of $V$ and hence also of $Y$.
This distribution and the contaminated normal family have been employed by
\cite{lang:sins:1993}
as the stochastic constituents for an adaptive robust regression methods.


\section{Variance-mean mixtures and their generalization}  \label{s:var-mean-mix}

The more versatile formulations are those where both
variables $R$ and $S$ in \ref{e:gmn-RS}, are non-degenerate, 
which constitute the theme of the present section.
 
\paragraph{Variance-mean normal mixtures}
The archetypal construction of  this form is the class  of distributions
with representation \ref{e:var-mean-mix}. In principle, any choice for
the distribution of the mixing variable $V$ is feasible, provided $\pr{V>0}=1$.

In practice, the predominant formulation occurs when $V$ follows a GIG distribution
with density \ref{e:pdf-gig}, leading to the GH family of distribution for $Y$, 
introduced by \citet{barndorff-nielsen:1977prs,barndorff-nielsen:1978sjs}.
The GH distributions depends of the set parameters
$\theta=(\lambda, \chi, \psi, \xi,\Sigma, \gamma)$,
and its density function at  $y\in\Real^d$ is%
\begin{eqnarray}
f_{\GH}(y; \theta)
   = \:\frac{(\sqrt{\psi/\chi})^\lambda \:C^{\frac{d}{2}-\lambda}}%
      {(2\pi)^{d/2}\det(\Sigma)^{1/2}K_\lambda(\sqrt{\chi\psi})}
 \:
 \frac{K_{\lambda-\frac{d}{2}}\left(\sqrt{C\:(\chi+\|y-\xi\|^2_\Sigma)}\right)
  \:e^{\gamma\T \Sigma\inv (y-\xi)}}%
  {\left(\sqrt{C\:(\chi+\|y-\xi\|^2_\Sigma)}\right)^{\frac{d}{2}-\lambda}}
  \label{e:pdf-gh}
\end{eqnarray}
where we have used the notation \ref{e:M-dist} and set
$ C = \psi+ \|\gamma\|^2_\Sigma = \psi+\alpha^2$,
bearing in mind \ref{e:aux-param}.

As presented, the parametric family is non-identifiable, due to the coexistence
of a scale factor in the GIG distribution and of a similar component in $\Sigma$.
The problem is resolved by imposing some restriction on the parameters.
A classical choice, adopted in the original Barndorff-Nielsen's papers,
is to set $\det(\Sigma)=1$.

After the initial formulation mentioned above, the GH family has been developed
further in a series of papers, and it has been employed in much applied work
across a wide range of domains.
Given the popularity of the GH family, we do not attempt to discuss its wide
body of properties and uses.
However, we must at least recall here that the GH family is
closed with respect to a number of operations, namely  marginalization,
conditioning and affine transformations; see \cite{blaesild:1981}.
For additional information, we refer to the detailed overview of this
and related distributions  presented in Chapter~9 of \cite{paolella:2007}.
Another account, inclusive of numerical illustrations in quantitative finance,
is available in Section~3.2 of \cite{mnei:frey:embr:2005};
we note  their mention on p.\,79 of the ``bewildering array of
alternative parameterizations''  introduced in the literature.

In our formulation,  \ref{e:var-mean-mix} corresponds to the specification
$r(u,v)=v$, $s(u,v)=v^{1/2}$ in \ref{e:gmn}, or equivalently $R=V$,
$S=V^{1/2}$ in \ref{e:gmn-RS}, with only one mixing variable involved.


Although the choice of a GIG distribution for $V$ in \ref{e:var-mean-mix}
is the predominant one, this does not rule out other possibilities.
In the univariate context, \cite{sichel:1973} had proposed the
formulation \ref{e:var-mean-mix} but assuming a Gamma distribution for $V$.
Given the above-mentioned fact that the Gamma family is a boundary case
of the GIG family, a multivariate version of the original Sichel's
construction can be obtained by setting $\chi=0$ and $\lambda>0$
in \ref{e:pdf-gh}.

\paragraph{Scale mixtures of skew-normal distributions}
\citet[Section 3.1]{bran:dey:2001} have introduced, with a very slightly
different name,   the idea of scale mixtures of SN variables.
These can be represented as $Y=\xi+ V^{1/2}Z$, where $Z$ has density
of type \ref{e:pdf-esn} with $\tau=0$ and $\xi=0$.
Combining this representation with the additive representation \ref{e:mean-mix}
of skew-normal variables, namely $Z= U\gamma + X$  where $U$ is half-normal,
we write
\begin{equation}
    Y =\xi+ V^{1/2}Z= \xi + V^{1/2}\:(U\gamma + X)
      = \xi + U\,V^{1/2} \gamma + V^{1/2} X
    \label{e:sn-scale-mix}
\end{equation}
which is of type \ref{e:gmn} with $r(u,v)= u\,v^{1/2}$ and $s(u,v)=v^{1/2}$.

A member of this class which has received much attention since 2001,
both on  the theoretical and the applied side, is the skew-$t$ (ST)
distribution which occurs when $V\sim \nu/\chi^2_\nu$.
A member of the ST family is therefore identified by the set of
parameters $(\xi, \Sigma, \gamma, \nu)$. The expression of the ST
density can be obtained as a special case of density \ref{e:pdf-est}
below, when $\tau=0$.
Chapters 4 and~6 of \citet{azza:capi:2014} provide a fairly detailed
account of this distribution, in the univariate and the multivariate case.
Note that a different parameterization is adopted there.

Another instance of scale mixtures of SN variates is represented
by the skew-slash distribution proposed by \cite{wangJ:gent:2006}.
Similarly to the symmetric multivariate slash distribution,
here $V=W^{-r}$ where $W\sim U(0,1)$ and $r$ is a positive parameter.

An instance of \ref{e:sn-scale-mix} producing to a very broad
set of distributions has been examined by \cite{vilc:bala:zell:2014}
taking $V$ to be a GIG variable, leading to what they denote
`skew-normal generalized hyperbolic distribution'.
Analogously to the classical GH  distribution \ref{e:pdf-gh},
this family is extremely flexible and it includes as a special
case several existing parametric families, as illustrated
in detail in Section~3 of the quoted paper.
A suitable constraint is required to avoid overparameterization,
such as $\det(\Sigma)=1$ employed for the classical GH distribution.

\paragraph{Extended skew-$t$ distribution} As the name indicates, the
extended ST family is a superset of ST family, introduced in independent
work by \cite{adcock:2010} and by \cite{arel:gent:2010metron}.
Despite the close connection with the ST family and many common properties,
the stochastic representation of extended version in form \ref{e:gmn}
involves a different mechanism from the one of the ST distribution.

Introduce a random variable $U$ having univariate $t$  distribution
on $\nu$ degrees of freedom truncated below $-\tau$, where $\tau\in\Real$
represents an additional parameter.
Also, let $V\inv\sim \chi^2_\nu/\nu$ like for the ST.
Then, taking into account Proposition~2 of \cite{arel:gent:2010metron},
we can say that
\begin{equation}
  Y = \xi + U\gamma + \left(\frac{\nu + U^2}{\nu+1} \:V\right)^{1/2} X
  \label{e:est}
\end{equation}
has extended ST distribution with density at $y\in\Real^d$, 
regulated by parameters $(\xi, \Sigma, \gamma, \nu, \tau)$,  equal to
\begin{equation}
\frac{1}{T(\tau;\nu)}t_d\left(y;\xi,\Omega,\nu\right)
  T\left(\left(\bar\tau+\eta\T(y-\xi)\right)
  \left(\frac{\nu+1}{\nu+\|y-\xi\|^2_\Omega}\right)^{1/2};\nu+1\right)
 \label{e:pdf-est}
\end{equation}
where $\Omega$, $\eta$ and $\bar\tau$ are as in
\ref{e:aux-param}, and $T(x;\nu)$ denotes the distribution function of
a standard univariate Student's $t$ variable on $\nu$ degrees of freedom.

\section{Final remarks} \label{s:final}

As anticipated in Section~\ref{s:normal-mix}, the main aims of the present
work are: (i)~to present a wide formulation, denoted GMN, which encompasses 
a large number of existing constructions involving continuous mixtures of normal
variables,  possibly in an implicit way; 
(ii)~to show that a unifying treatment is possible, by providing a number 
of general properties for the GMN class.

It would be possible to extend further this construction, even considerably.
Among the various options, a simple one would be along the following lines.
Start from the familiar construction represented
by the class \ref{e:var-mix} of scale mixtures of normals, and recall that 
a vast subset of the elliptical class of distributions can be expressed 
using scale mixtures of normal variates. 
Then combine this mechanism with an additive term like $V\,\gamma$ in 
\ref{e:var-mean-mix} or, similarly, $R\,\gamma$ in \ref{e:gmn-RS}.
Not only this extension would be possible, but a number of properties
developed in Section~\ref{s:gmn} would even carry on, for instance those 
in Subsections~\ref{s:quadratic-forms}  and~\ref{s:beta-Mardia}.%

There are, however, other facts which would not be preserved in this
extended construction. We are referring specifically to the properties of
closure under affine transformations and marginalization examined
in  Subsection~\ref{s:distrib-properties}. 
To see the source of the problem, consider a specific but fairly 
popular instance, namely the so-called exponential power distribution,
but also other names are in use.
The univariate formulation of \cite{subbotin:1923}  has subsequently been 
extended to the multivariate setting and this distribution can be 
represented as a scale mixture of normal variables, as for the parameter
set which corresponds to leptokurtic distributions, with heavier-than-normal 
tail behaviour. An explicit expression of the implied mixing distribution
is given by \cite{gome:gome:marin:2008}; a crucial fact is that this mixing
distribution depends of the dimension $d$ of the mixed normal distribution, 
$d$ in our notation. This situation prevents the property of closure 
under marginalization for a scale mixture of normal variables, 
as shown in the already-quoted work of \cite{kano:1994}, who has focused 
precisely on the case of the exponential power distribution.
Therefore, \emph{a fortiori}, closure of the class under the more
general manipulation of affine transformations cannot hold, even less so
if an additive stochastic term $U\gamma$ is included like in \ref{e:gmn-RS}.

To summarize, an extension of the construction along the above-delineated 
lines would be possible, and it would even preserve certain formal properties.
Other properties would not carry on, especially closure under affine 
transformations and marginalization, and these seem to us important facts
when we come the use of these constructions in applied work.
Obviously, this statement cannot be taken as definitive and an
absolute bar, since every such judgement must by evaluated with respect
to a given context. 
However, the above discussion explains why, in all the special cases
which we have examined, neither the mixing distribution of $(U, V)$ 
nor the functions $r$ and $s$ in \ref{e:gmn} depend on $d$.

Another possible direction for extension of the GMN construction is via 
consideration of more than two mixing variables. 
To exemplify in a simple form, one option would be to replace the 
term $U\gamma$ in \ref{e:mean-mix} by $\Gamma\,U$ where $\Gamma$ is
a matrix of coefficients and $U$ is a random vector.
This step would allow to incorporate families such as the so-called 
closed/unified skew-normal discussed by \cite{arel:azza:2006}.
Similarly, the scale mixture construction \ref{e:var-mix} can be 
extended by consideration of multiple random scale factors,
which amounts to the `multiple scale mixtures' proposed by
\cite{forb:wrai:2014}.
Combination of these two mechanisms would extend \ref{e:gmn-RS} to
the form
\[
   Y = \xi + \Gamma\,R + S\, X
\]
where $S$ is a diagonal matrix formed by $d$ positive random
variables. Clearly, such a study represents a separate undertaking,
with non-trivial complexity if one attempts to develop a unified 
treatment of the connected properties.

As for the adopted term `generalized mixtures of normal variables',
one could perhaps object that the term may suggests more than it actually means, 
since the construction does not encompass all possible mixtures of normal variables. 
The term must rather be intended in the same spirit of other similar instances, 
such as `generalized linear models', which do not embrace all possible 
extensions of linear models.



\appendix
\section*{Appendix}

\paragraph{Derivation of density function \ref{e:pdf-esn}}
By assumption, $(Y|U=u)\sim \N_d(\xi+\gamma u,\Sigma)$, and the density function of
$U$ is $2\phi(u)I_{(0,\infty)}(u)$.
With these positions, density \ref{e:pdf-gmn} becomes
\begin{eqnarray*}
f(y;\xi,\gamma,\Sigma)&=&\int_{0}^{\infty}\phi_d(y;\xi+u\gamma,\Sigma)\:2\:\phi(u)\d{u}\\
&=& 2\:\phi_d\left(y;\xi,\Sigma+\gamma\gamma\T\right)\:
 \int_{0}^{\infty}\phi\left(u;\gamma\T(\Sigma+\gamma\gamma\T)\inv(y-\xi),1-\gamma\T(\Sigma+\gamma\gamma\T)\inv\gamma\right) \d{u}\\
&=& 2\:\phi_d\left(y;\xi,\Sigma+\gamma\gamma\T\right)\int_{0}^{\infty}\:
 \phi\left(u;\frac{\gamma\T\Sigma\inv(y-\xi)}{1+\gamma\T\Sigma\inv\gamma},\frac{1}{1+       \gamma\T\Sigma\inv\gamma}\right) \d{u}\,.
 \end{eqnarray*}
taking into account the follow identity:
\[
\phi_d(y;\xi+u\gamma,\Sigma)\:\phi_1(u) \; = \;
\phi_d\left(y;\xi,\Sigma+\gamma\gamma\T\right)\:
\phi\left(u;\gamma\T(\Sigma+\gamma\gamma\T)\inv(y-\xi), 1-\gamma\T(\Sigma+\gamma\gamma\T)\inv\gamma\right) \,,
\]
and the inverse of $\Sigma + \gamma\gamma\T$ has been expressed
using the Sharman-Morrison formula \ref{e:SM}.
After the change variable
\[
v=\sqrt{1+\gamma\T\Sigma\inv\gamma}\left(u-\frac{\gamma\T\Sigma\inv(y-\xi)}{1+\gamma\T\Sigma\inv\gamma}\right)
\]
we can re-write
\[
f(y;\xi,\Sigma,\gamma) = 2\:\phi_d\left(y;\xi,\Sigma+\gamma\gamma\T\right)
  \int_{-\frac{\gamma\T\Sigma\inv(y-\xi)}{\sqrt{1+\gamma\T\Sigma\inv\gamma}}}^{\infty}\phi(v) \d{v}
\]
which coincides with density \ref{e:pdf-sn} under the notation defined in \ref{e:aux-param}.


\paragraph{Conditional distribution of multivariate Student's $t$ components}

Given a partition of $Y$ in  sub-vectors $Y_1$ and $Y_2$ as in \ref{e:partition},
we want to apply the general expressions \ref{e:hc}--\ref{e:f1|2} to find
the conditional distribution of $Y_1$ given that $Y_2{=}y_2$ for a vector $y_2\in\Real^{d_2}$,
in the special case when $Y$ has a $d$-dimensional Student's $t$ density \ref{e:pdf-t}.

With respect to the general expression \ref{e:gmn-RS}, here $R\equiv0$,
$S=V^{1/2}$ with $V\sim\nu/\chi^2_\nu$. 
It is immediate that the density of $V$ is
\[
  h(v)=\frac{(\nu/2)^{\nu/2}}{\Gamma(\nu/2)}\,v^{-(\nu+2)/2}e^{-\nu/(2v)}, \qquad v>0.
\]
which here plays the role of $h(r,s)$ in \ref{e:hc}.
Since $(Y_2|V=v)\sim \N_{d_2}(\xi_2,v\Sigma_{22})$,
it follows that $Y_2=V^{1/2}X_2\sim t_{d_2}(\xi_2,\Sigma_{22},\nu)$, so that we write
\begin{eqnarray*}
 \phi_{d_2}(y_2;\xi_2+v\Sigma_{22}) &=&
   \det(\Sigma_{22})^{-1/2} (2\pi)^{-d_2/2}v^{-d_2/2} e^{-q_2(y_2)/(2v)},\\
f_2(y_2)
&=&\frac{\det(\Sigma_{22})^{-1/2}\Gamma((\nu+d_2)/2)\nu^{\nu/2}}%
   {\Gamma(\nu/2)\pi^{d_2/2}}(\nu+q_2(y_2))^{-(\nu+d_2)/2}
\end{eqnarray*}
having set
\[
  q_2(y_2)= (y_2-\xi_2)\T\Sigma_{22}\inv(y_2-\xi_2) = \|y_2-\xi_2\|^2_{\Sigma_{22}}\,.
\]
After plugging these ingredients in \ref{e:hc} and some simplification, we obtain  that
\begin{eqnarray*}
  h_c(v|y_2)
  &=&\frac{[(\nu+q_2(y_2))/2]^{(\nu+d_2)/2}v^{-(\nu+d_2+2)/2} e^{- (\nu+q_2(y_2))/(2v)}}{\Gamma((\nu+d_2)/2)}\:,
  \qquad v>0,
\end{eqnarray*}
which means that $(V|Y_2=y_2)\sim(\nu+q_2(y_2))/\chi^2_{\nu+d_2}$.
This distribution  is the same of
$(\nu+d_2)\inv(\nu+q_2(y_2))\, V_2$,
where $V_2\sim (\nu+d_2)/\chi^2_{\nu+d_2}$.
Hence, by setting $X_{1|2}\sim \N_{d_1}(0,\Sigma_{11|2})$ be  a variable independent of $V_2$,
we obtain
\[
   (Y_1|Y_2=y_2)   \equald 
      \xi_{1|2} +\left(\frac{\nu+q_2(y_2)}{\nu+d_2}\right)^{1/2}\,V_2^{1/2} X_{1|2}.
\]
Since $V_2^{1/2} X_{1|2}\sim t_{d_1}(0,\Sigma_{11|2},\nu+d_2)$,  we conclude that
\[
  (Y_1|Y_2=y_2)  \sim  
  t_{d_1}\left(\xi_{1|2},\left(\frac{\nu+q_2(y_2)}{\nu+d_2}\right)\Sigma_{11|2},\nu+d_2\right)
\]
indicating that $(Y_1|Y_2=y_2)$ is still of Student's $t$ type.

If we insert the parameters of this distribution in the density \ref{e:pdf-t},
we obtain an expression algebraically equivalent to formula (1.15) of \cite{kotz:nada:2004}
for the conditional $t$ density.
It could possibly be remarked that our derivation is not any simpler than the one of  \cite{kotz:nada:2004},
and perhaps even a little more lengthy, but it has the advantage of revealing the
$t$ nature of this density, a fact not so visible from the other development.

\paragraph{Derivation of density function \ref{e:pdf-Rayleigh-mix}}

The conditional density of $(Z\mid U = u)$ is $\phi_d(z;\gamma u,\Sigma)$.
Hence, its unconditional density  is
\begin{eqnarray*}
 f_Z(z;\gamma,\Sigma) &=& \int_{0}^{\infty} \phi_d(z;\gamma u,\Sigma) g_u(u) \d u \\
 &=& \int_{0}^{\infty} (2\pi)^{-d/2} \det(\Sigma)^{-1/2} e^{-\half (z-\gamma u)\T\Sigma\inv(z-\gamma u)} u e^{-\half u^2} \d u \\
 &=& (2\pi)^{-d/2}\det(\Sigma)^{-1/2}\int_{0}^{\infty} u\,e^{-\half \{u^2 + (z-\gamma u)\T\Sigma\inv(z-\gamma u)\}} \d u.
\end{eqnarray*}
By making use of the identities,
\begin{eqnarray*}
& u^2 + (z-\gamma u)\T\Sigma\inv(z-\gamma u)
= z\T\Omega\inv z + (1 + \alpha^2)(u-\tilde\eta\T z)^2\,,& \\
&\det(\Omega) = \det(\Sigma+\gamma\gamma\T)= \det(\Sigma)\:(1+\alpha^2)\,,&
\end{eqnarray*}
where $\Omega$, $\eta$ and $\alpha$ are as in \ref{e:aux-param}, and $\tilde\eta=(1+\alpha^2)^{-1/2}\eta$,
we have
\begin{eqnarray*}
 f_Z(z;\gamma,\Sigma) &=& (2\pi)^{-d/2}\det(\Sigma)^{-1/2}e^{-\half z\T\Omega\inv z}\int_{0}^{\infty} u\,e^{-\half (1 + \alpha^2)(u-\tilde\eta\T z)^2} \d u\\
 &=& (2\pi)^{-d/2} \det(\Omega)^{-1/2}e^{-\half z\T\Omega\inv z}\int_{0}^{\infty} u (1+\alpha^2)^{1/2} e^{-\half (1 + \alpha^2)(u-\tilde\eta\T z)^2} \d u\\
 &=& (2\pi)^{1/2}\phi_d(z;\Omega) \int_{0}^{\infty} u (1+\alpha^2)^{1/2} \phi((1+\alpha^2)^{1/2} (u-\tilde\eta\T z))\d u.
\end{eqnarray*}
The change variable $w=(1+\alpha^2)^{1/2} (u-\tilde\eta\T z)$ yields
\begin{eqnarray*}
 f_Z(z;\gamma,\Sigma)  &=& (2\pi)^{1/2}\phi_d(z;\Omega) \int_{-\eta\T z}^{\infty}
 \left\{(1+\alpha^2)^{-1/2} w + \tilde\eta\T z\right\} \phi(w)\d w\\
 &=& (2\pi)^{1/2}(1+\alpha^2)^{-1/2}\phi_d(z;\Omega)  \left\{\int_{-\eta\T z}^{\infty} w \phi(w) \d w + \eta\T z
 \int_{-\eta\T z}^{\infty} \phi(w)\d w \right\}\\
 &=& (2\pi)^{1/2}(1+\alpha^2)^{-1/2}\phi_d(z;\Omega) \left\{-\int_{-\eta\T z}^{\infty} \d\phi(w)  + \eta\T z
 \int_{-\infty}^{\eta\T z} \phi(w)\d w\right\}\\
 &=& (2\pi)^{1/2}(1+\alpha^2)^{-1/2}\phi_d(z;\Omega) \:\left\{\phi(\eta\T z)  + \eta\T z \Phi(\eta\T z)\right\}
 \end{eqnarray*}
which coincides with expression \ref{e:pdf-Rayleigh-mix}.


\paragraph{Derivation of the Mardia's multivariate measures of skewness and  kurtosis}

We expand here the computations sketched in  Subsection~\ref{s:beta-Mardia} for
computing the Mardia's measures.
Given the positions stated in the initial part of Subsection~\ref{s:beta-Mardia},
we can rewrite
\[\beta_{d,1}=\E{[Y_0\T\Sigma^{1/2}\Sigma_Y\inv \Sigma^{1/2}Y_0']^3}\qquad\beta_{d,1}=\E{[Y_0\T\Sigma^{1/2}\Sigma_Y\inv\Sigma^{1/2} Y_0]^2},\]
where $Y_0'=\Sigma^{-1/2}(Y'-\mu_Y)$ and we can expand
\[
\Sigma^{1/2}\Sigma_Y\inv \Sigma^{1/2}
= \mu_{02}\inv\Sigma^{1/2}
 \left(\Sigma\inv-\frac{\rho}{1+\rho\alpha^2}\Sigma\inv\gamma\gamma\T\Sigma\inv\right)
\Sigma^{1/2}
= \mu_{02}\inv\left(I_d-\bar\rho\gamma_0\gamma_0\T\right)
\]
using the Sherman-Morrison formula \ref{e:SM} to invert $\Sigma_Y$
and denoting $\mu_{02}=\E{S^2}$. Hence, re-write further
\[\beta_{d,1}=\mu_{02}^{-3}\E{[Y_0\T(I_d- \bar\rho\gamma_0\gamma_0\T)Y_0']^3}
=\mu_{02}^{-3}\E{[Y_0\T Y_0'-\bar\rho (\gamma_0\T Y_0)(\gamma_0\T Y_0')]^3},\]
and
\[\beta_{d,2}=\mu_{02}^{-2}\E{[Y_0\T(I_d-\bar\rho\gamma_0\gamma_0\T)Y_0]^2}
=\mu_{02}^{-2}\E{[\|Y_0\|^2-\bar\rho(\gamma_0\T Y_0)^2]^2}.\]
On setting $T_0=\alpha\inv W_0$, where $W_0=\gamma\T\Sigma\inv X$, we note that
\[
  \gamma_0\T Y_0
    =\gamma\T\Sigma\inv(R_0\gamma+S\,X)
    =R_0 \alpha^2 +S\,\gamma\T\Sigma\inv X=\alpha(R_0\alpha+S\,T_0)
    =\alpha Z_0,
 \]
where $Z_0=R_0\alpha+S\, T_0$, with  $(R_0,S)$ and $T_0$ independent variables.
Thus, by letting  $Z_0'$ be an  independent copy of $Z_0$,
we have that
\[\beta_{d,1}=\mu_{02}^{-3}\E{[Y_0\T Y_0'-\bar\rho Z_0 Z_0')]^3}=\mu_{02}^{-3}\E{[Y_0\T Y_0'-Z_0 Z_0'+(1-\bar\rho) Z_0 Z_0')]^3},\]
and
\[\beta_{d,2}=\mu_{02}^{-2}\E{[\|Y_0\|^2-\bar\rho Z_0^2]^2}=\mu_{02}^{-2}\E{[\|Y_0\|^2-Z_0^2+(1-\bar\rho) Z_0^2]^2}.\]
Note that 
\[Y_0\T Y_0'=(R_0\gamma_0+S\,X_0)\T(R_0'\gamma_0+ S' \,X_0')
  =\alpha^2 R_0R_0'+\alpha(R_0 S' \,T_0'+R_0'S\,T_0)+ S\,S' \,X_0\T X_0'
\]
and
\[Z_0 Z_0'=(R_0\alpha+S\, T_0)(R_0'\alpha+ S' \, T_0')=\alpha^2 R_0R_0'+\alpha(R_0 S' \,T_0'+R_0'S\,T_0)+ S\,S' \,T_0T_0'\]
leading to
\[Y_0\T Y_0'-Z_0 Z_0'= S\,S' \,(X_0\T X_0'-T_0 T_0')= S\,S' \,X_0\T M_0 X_0',\]
where $M_0=I_d-\bar\gamma\bar\gamma\T$ is a projection matrix.
Similarly, we find
\[\|Y_0\|^2-Z_0^2=S(\|X_0\|^2-T_0^2)= S^2\, X_0\T M_0X_0.\]

Recall that $M_0X_0$  and $T_0$  are independent, and so are  $M_0X_0'$ and $T_0'$.
Therefore, $M_0X_0$ is also independent of $Z_0$ since this variable depends on $(R_0, S, T_0)$ only;
similarly, independence holds for $M_0X_0'$ and $Z_0'$.
Finally, we find that
\[\beta_{d,1}=\mu_{02}^{-3}\E{[ S\,S' \,(X_0\T M_0 X_0')+(1-\bar\rho) Z_0 Z_0')]^3},\]
and
\[\beta_{d,2}=\mu_{02}^{-2}\E{[S^2\,(X_0\T M_0X_0)+(1-\bar\rho) Z_0^2]^2},\]
where we note that $MX_0$ and $Z_0$ are independent random quantities of mean zero.

We can now start calculation of the Mardia's measures:
\begin{eqnarray*}
\beta_{d,1}&=&\mu_{02}^{-3}\E{[ S\,S' \,(X_0\T M_0 X_0')+(1-\bar\rho) Z_0 Z_0')]^3}\\
&=&\mu_{02}^{-3}[\E{ S\,S' \,(X_0\T M_0 X_0')^3}+3(1-\bar\rho)\E{ S\,S' (X_0\T M_0 X_0')^2(Z_0 Z_0')}\\
&& \qquad +3(1-\bar\rho)^2\E{ S\,S' \,(X_0\T M_0 X_0')(Z_0 Z_0')^2}+(1-\bar\rho)^3
\E{(Z_0 Z_0')^3}]\\
&=&\mu_{02}^{-3}[\E{ S\,S' \,}\E{(X_0\T M_0 X_0')^3}+3(1-\bar\rho)\E{(X_0\T M_0 X_0')^2}\E{ S\,S' (Z_0 Z_0')}\\
&& \qquad +3(1-\bar\rho)^2\E{(X_0\T M_0 X_0')}\E{ S\,S' \,(Z_0 Z_0')^2}+(1-\bar\rho)^3
\E{(Z_0 Z_0')^3}],
\end{eqnarray*}
where by symmetry $\E{(X_0\T M_0 X_0')^3}=\E{(X_0\T M_0 X_0')}=0$,  and
\begin{eqnarray*}
\E{(X_0\T M_0 X_0')^2}&=&\E{\tr\{X_0\T M_0 X_0'(X_0')\T M_0 X_0\}}\\
&=&\tr\{\E{M_0 X_0'(X_0')\T M_0 X_0X_0\T}\} \\
&=&\tr\{M_0\E{X_0'(X_0')\T} M_0\E{X_0X_0\T}\} \\
&=&\tr(M_0^2)=\tr(M_0)=d-1.
\end{eqnarray*}
where we have use the fact $\E{X_0'\,(X_0')\T}=\E{X_0X_0\T}=I_d$.
Thus, since $(S',Z_0')$ and $(S,Z_0)$ are independent variables with the same distribution,
 $\beta_{1,d}$ reduces to
\[
  \beta_{1,d}=\mu_{02}^{-3}
        [3(d-1)(1-\bar\rho)(\E{S^2\,Z_0})^2 + (1-\bar\rho)^3(\E{Z_0^3})^2],
\]
and analogously
\begin{eqnarray*}
\beta_{2,d}
&=&\mu_{02}^{-2}\E{[S^2 (X_0\T M_0X_0)+(1-\bar\rho) Z_0^2]^2}\\
&=&\mu_{02}^{-2}[\E{S^4 (X_0\T M_0X_0)^2+2(1-\bar\rho)S^2 (X_0\T M_0X_0)Z_0^2+(1-\bar\rho)^2 Z_0^4}]\\
&=&\mu_{02}^{-2}[\E{S^4}\E{(X_0\T M_0X_0)^2}+2(1-\bar\rho)\E{(X_0\T M_0X_0)}\E{S^2\,Z_0^2}+(1-\bar\rho)^2\E{Z_0^4}]\\
&=&\mu_{02}^{-2}[(d+1)(d-1)\E{S^4}+2(d-1)(1-\bar\rho)\E{S^2\,Z_0^2}+(1-\bar\rho)^2\E{Z_0^4}]
\end{eqnarray*}
considering that $X_0\T M_0X_0\sim\chi_{d-1}^2$.

To complete the calculations, we need some moments of the distribution of $(S,Z_0)$.
For this, recall that  $(R_0, S)$ and $T_0\sim \N(0,1)$ are independent
variables, and (i)~the odd-order moments of $T_0$ are $0$,
(ii)~$\E{R_0}=0$,  $\E{T_0^2}=1$, $\E{T_0^4}=3$, and
(iii)~we are assuming that both $R$ and $S$ have finite moments
up to order four. Then, we need just the following expected values:
\begin{eqnarray*}
\E{S^2\,Z_0}
&=&\E{\alpha R_0\,S^2+S^3\,T_0}\\
&=&\alpha\E{R_0\,S^2}\\
&=&\alpha(\E{R\,S^2}-\E{S^2}\E{R}),\\
\E{S^2\,Z_0^2}&=&\E{\alpha^2 R_0^2\,S^2+2\alpha R_0 S^3\,T_0+S^4\,T_0^2}\\
&=&\alpha^2\E{R_0^2\,S^2}+\E{S^4}\\
&=&\alpha^2(\E{R^2\,S^2}-2\E{R\,S^2}\E{R}+(\E{R})^2\E{S^2})+\E{S^4},\\
\E{Z_0^3}
&=&\E{\alpha^3R_0^3+3\alpha^2R_0^2 S\,T_0+3\alpha R_0 S^2 T_0^2+S^3 T_0^3}\\
&=&\alpha^3\E{R_0^3} + 3\alpha\E{R_0 S^2}\\
&=&\alpha^3(\E{R^3}-3\E{R^2}\E{R}+2(\E{R})^3) + 3\alpha(\E{R\,S^2}-\E{R}\E{S^2}),\\
\E{Z_0^4}
&=&\E{\alpha^4R_0^4+4\alpha^3R_0^3S\,T_0+6\alpha^2R_0^2ST_0^2+4\alpha R_0 S^3\,T_0^3+S^4\,T_0^4}\\
&=&\alpha^4\E{R_0^4}+6\alpha^2\E{R_0^2\,S^2}+3\E{S^4}\\
&=&\alpha^4(\E{R^4}-4\E{R^3}\E{R}+6\E{R^2}(\E{R})^2-3(\E{R})^4)\\
&&+6\alpha^2(\E{R^2\,S^2}-2\E{R\,S^2}\E{R}+(\E{R})^2\E{S^2})+3\E{S^4}.
\end{eqnarray*}
After substitution of various expectations with the symbols defined in \ref{e:mu.hk},
we arrive at the expressions reported in Subsection~\ref{s:beta-Mardia}.


\end{document}